\documentclass[11pt]{article}

\usepackage[pdftex]{graphicx}
\usepackage[ddmmyyyy]{datetime}
\usepackage{amssymb,amsfonts,amsmath}
\usepackage{amsthm,xypic,enumerate}
\usepackage{color,tikz,pgf,subfigure}
\usepackage{fancyhdr}
\usetikzlibrary{shapes,arrows} 

\input xy
\xyoption{all}
\usepackage{hyperref}

\newenvironment{enumerate*}[1][{}]{\begin{itemize}\renewcommand{\labelitemi}{#1}}{\end{itemize}}

\numberwithin{equation}{section}
\newtheorem{theorem}[equation]{Theorem}
\newtheorem{proposition}[equation]{Proposition}
\newtheorem{corollary}[equation]{Corollary}
\newtheorem{lemma}[equation]{Lemma}

\theoremstyle{definition}
\newtheorem{definition}[equation]{Definition}
\newtheorem{remark}[equation]{Remark}
\newtheorem{example}[equation]{Example}

\renewenvironment{proof}[1][{}]{\medskip\noindent{\bf Proof of {#1}.}}{\qed}

\newcommand{\N}{\mathbb{N}}
\newcommand{\R}{\mathbb{R}}
\newcommand{\mmS}{\mathcal{S}}
\newcommand{\mmC}{\mathcal{C}}
\newcommand{\mmN}{\mathcal{N}}
\newcommand{\mmR}{\mathcal{R}}
\newcommand{\omR}{\overline{\R}}

\newcommand{\omN}{\overline{\N}}
\newcommand{\wm}{\widetilde{M}}
\newcommand{\mmK}{\mathcal{K}}

\DeclareMathOperator*{\im}{im}

\DeclareMathOperator*{\supp}{supp}
\DeclareMathOperator*{\sign}{sign}
\DeclareMathOperator*{\diag}{diag}

\newcommand{\ta}{\tilde{a}}
\newcommand{\tb}{\tilde{b}}

\title{Power-law Kinetics and Determinant Criteria for the Preclusion of Multistationarity in Networks of Interacting Species}

\author{Carsten Wiuf\footnotemark[2] \and Elisenda Feliu\footnotemark[2]}

\begin{document}

\maketitle

\renewcommand{\thefootnote}{\fnsymbol{footnote}}
\footnotetext[2]{Department of Mathematical Sciences, University of Copenhagen, Universitetsparken 5, 2100 Denmark.}
\renewcommand{\thefootnote}{\arabic{footnote}}

\begin{abstract} 
We present determinant criteria for the preclusion of non-degenerate multiple steady states in  networks of interacting species. 
A network is modeled as a system of ordinary differential equations in which the form of the species formation  rate function is restricted by the reactions of the network and how the species influence each reaction.
We characterize families of so-called power-law kinetics for which the associated species formation rate function is injective within each stoichiometric class and thus the  network cannot exhibit multistationarity.
The criterion for power-law kinetics is derived from the determinant of the Jacobian of the species formation rate function.
Using this characterization we further derive similar determinant criteria applicable to general sets of kinetics.
The criteria are conceptually simple, computationally tractable  and easily implemented.
 Our approach embraces and extends  previous work on multistationarity, such as work in relation to chemical reaction networks with dynamics defined by  mass-action or non-catalytic  kinetics, and also work based on  graphical analysis of the interaction graph associated to the system. Further, we interpret the criteria in terms of circuits in the so-called DSR-graph.

\end{abstract}

\pagestyle{myheadings}
\thispagestyle{plain}
\markboth{C. Wiuf, E. Feliu}{Preclusion of Multiple Steady States}

\section{Introduction}

Networks of interacting species are used in many areas of science to represent the structural form of a dynamical system. This is in particular the case in systems biology and biochemistry where biochemical reactions are represented in the form of a network. However, similar network structures are also used  in  ecology, cell biology and epidemics, as well as outside the natural sciences, to describe the possible interactions between some species of interest. Common to these  networks is that they consist of a set of species and a set of interactions among the species. The state of the system is given by the concentration (or abundance) of each species and each interaction represents a transformation of the state of the system.  An example is the chemical reaction $A+B\rightarrow 2C$ where one molecule of $A$ and one molecule of $B$ form two molecules of $C$.

Typically, a  system of ordinary differential equations (ODEs) is used to describe how  species concentrations change over time. The species formation rate function of the system  describes  the instantaneous change in the concentrations when considering simultaneously the individual rates of all reactions in the network. Reaction rates are generally unknown but some qualitative aspects might be assumed or inferred. For instance, the presence of the species on the left side of a reaction ($A,B$ in the above example)  might be a prerequisite for the reaction to take place and higher concentrations of these species typically lead to higher reaction rates. 
In some cases, reaction rates are fixed to follow a specific functional form that might  depend on parameters to be inferred from experimental observations, in addition to the concentrations of the species. In other cases, only weak assumptions are imposed on the reaction rates and functions. We consider restrictions given by a so-called influence specification \cite{shinar-conc}, which specifies how each species affects (positively, negatively, or neutrally) the reaction rates with increasing concentration.
As pointed out by other authors \cite{banaji-craciun2}, the full structure of the network (in particular the reactions) is not required to develop the theory. We will, however,  keep the terminology of reaction networks, as they provide the main source of inspiration and examples.

For many networks the structure of the interactions and the influence specification alone determine dynamical and steady-state properties of the system (for example, multistationarity, persistence, or oscillations). That is to say, irrespectively the rates
and the parameters quantifying them, taken together with the initial species concentrations, the system shows qualitatively the same type of behavior. It is  perhaps surprising as the network structure itself does not encode any information about the specific rate functions and abundances. Even small networks might have many parameters  which potentially could give rise to a rich and varied dynamics, as well as differences in  the long-term behavior of the system.

Of particular interest has been to determine whether a system allows for multiple positive steady states, also known as multistationarity. Multistationarity provides a mechanism for switching (rapidly) between different  responses and confers robustness to the steady-state values of the system \cite{Huang-Ferrell,Markevich-mapk}. One way to address whether a system exhibits multistationarity is by finding the positive solutions to the steady-state equations of the system. 
Solving the equations might prove difficult, if not impossible, with   difficulty depending on the assumptions about the reactions rates and the number of parameters.
Here we take a more conceptual route and focus on understanding the characteristics of   networks that cannot exhibit multistationarity, irrespectively of the specific choices of reaction rates.

Various criteria have been developed  to preclude the existence of multiple positive steady states for general classes of reaction functions, also called  kinetics, such as mass-action kinetics \cite{Feinbergss,craciun-feinberg,craciun-feinbergI,Feliu-inj,feinberg-def0,conradi-switch,conradi-PNAS}, non-catalytic kinetics \cite{banaji-donnell}, and weakly monotonic kinetics \cite{shinar-conc}.  These criteria typically utilize the structure  of the system together with some assumptions about the form of the rate functions. For example, for mass-action kinetics the rate functions are polynomials and the steady-state equations become a system of polynomial equations.  Capitalizing on the polynomial form of the equations  has lead to specific conditions to preclude multistationarity \cite{craciun-feinberg,craciun-feinbergI,Feliu-inj,PerezMillan}.

The aim of this paper is to provide a computationally tractable determinant criterion for \emph{injectivity} of a network 
for different classes of kinetics.  Injectivity refers  to injectivity of the 
species formation rate function that governs the dynamical behavior of the system. If this  function is injective for the allowed kinetics then  the system does not have the capacity for multiple positive steady states.  The idea of injectivity  was introduced by Craciun and Feinberg \cite{craciun-feinbergI}, but it is also  underlying previous work on the preclusion of multistationarity \cite{soule}. We show that injectivity of a network is closely related to injectivity of a network taken with \emph{power-law kinetics} \cite{hornjackson,clarke-SNA}. Power-law kinetics generalize  mass-action kinetics and confer greater flexibility to the form of the  rate functions than mass-action kinetics. Savageau \cite{savageau} emphasizes the importance of power-law kinetics  in biochemistry but their importance outside biology for modeling purposes is also well documented \cite{anderson-may,stroud}. It can be argued that power-law kinetics provide approximations to kinetics in general \cite{savageau}, which  is also exemplified in our work.

Power-law kinetics share common features with mass-action kinetics and  parallel results can be derived for the two types of kinetics.  We derive necessary and sufficient \emph{determinant criteria} for a network to be injective over different classes of power-law kinetics. The determinant refers to the determinant of a \emph{modified} version of the species  formation rate function  (Definition~\ref{extended2} in this paper). We show that the determinant of the Jacobian of the modified function is non-zero for all  concentration vectors and  kinetics in one of the classes  if and only if the network is injective over the particular class. For power-law kinetics our results extend parallel results for mass-action kinetics \cite{Feliu-inj}.  We proceed to show that injectivity over the class of all power-law kinetics compatible with an influence specification is equivalent to injectivity over the class of all kinetics that respect the same influence specification.  Therefore, we obtain determinant criteria to preclude the existence of multiple steady states for all kinetics that respect the influence specification, \emph{independently of the specific functional form of the kinetics} (Theorem~\ref{equiv}). The criteria depend in part on the network structure and in part on the influence specification. 

The criteria are easily implemented using symbolic software packages, such as Mathematica, Maple or SAGE, and, thus, they are of practical use.  We present two equivalent criteria: one involving the computation of a symbolic determinant, and the other involving the computation of minors of numerical matrices. For moderately-sized  networks the criteria are computationally efficient, and the first is usually faster.  For larger networks memory restrictions might constrain the  computation of  the symbolic determinant and time restrictions might constrain the computation using the minors.
Being based on the computation of determinants, the criteria can be reformulated in terms of circuits in a variant of the DSR-graph \cite{banaji-craciun1}. This might allow for the development of visual approaches to injectivity, in the style of  \cite{banaji-craciun1,banaji-craciun2,craciun-feinbergII,Shinar-conc2}.

Our work extends and embraces previous determinant criteria developed for networks taken with mass-action kinetics \cite{craciun-feinberg,craciun-feinbergI,Feliu-inj,craciun-feinberg-semiopen}. Further, it closely relates to recent work by Shinar and Feinberg \cite{shinar-conc}, where  a characterization of injective networks is provided  for classes of kinetics defined by an influence specification (these networks are called concordant networks). Their definition of influence specification differs from the one adopted here, but it can be recovered in our context (Section~\ref{sec:shinar}). Specifically, their definition corresponds to consider the union of certain classes of kinetics, rather than the classes individually. Instead of being determinant based, the criteria in \cite{shinar-conc} are based on computing the signs of vectors in different sets and hence the two approaches  differ in nature. Further, the present work clarifies the  role played by power-law kinetics in deciding  injectivity.
 
In a series of papers \cite{banaji-donnell,banaji-craciun1,banaji-craciun2}, the authors study    injectivity of a certain class of kinetics (called non-catalytic kinetics) and  of dynamical systems more generally. These articles tackle injectivity of so-called open networks (which contain all reactions of the form $S\rightarrow 0$, where $S$ is a species).  We provide a  discussion in Section~\ref{sec:banaji} of the relationship between  our results and  those in  \cite{banaji-donnell,banaji-craciun1,banaji-craciun2}.
Finally, this work also relates to a criterion for multistationarity based on the interaction graph  given by Kaufman, Soul\'{e} and Thomas \cite{kaufman,soule}. The interaction graph records the sign of the entries in the Jacobian of a dynamical system. In Section~\ref{sec:soule}, we relate our determinant criterion to that of \cite{kaufman,soule}.

The structure of the paper is the following. In Section \ref{sec:notation} we introduce some notation and in Section \ref{crns} we introduce the basic concepts of networks and kinetics  and Section \ref{sec:factor} presents the general form of the dynamical systems we consider. Section \ref{matrixresults} is concerned with some useful theoretical matrix results,  and Section \ref{sec:inj} introduces the notion of degeneracy and injectivity, two key concepts. Section \ref{gen-mass} and \ref{sec:powerlawinj} focus on power-law kinetics and derive a determinant criterion for injectivity.  In Section \ref{sec:int}, we discuss influence specifications, give examples from the literature and derive further results on injectivity for families of power-law kinetics. These results are extended in Section \ref{general-kinetics}  to broader and more general classes of kinetics. In Section \ref{graphical_sec} we develop a graph-theoretical interpretation of our criteria. Sections \ref{sec:shinar}, \ref{sec:banaji} and \ref{sec:soule} are devoted to  the relationship between our work and previous work \cite{shinar-conc,banaji-donnell,banaji-craciun2,kaufman}. Finally, in Section \ref{sec:hill} we show that other types of kinetics could be used in place of power-law kinetics. To keep the exposition clear in the main text, all proofs are in the Appendix.

\section{Notation}\label{sec:notation}

Let $\R_{+}$ denote the set of positive real numbers (without zero) and $\omR_{+}$  the set of non-negative real numbers (with zero). Similarly, let $\omN$ be the set of non-negative integers. 
Given a finite set $\mathcal{E}$,  the ring of polynomials in $\mathcal{E}$ is denoted $\R[\mathcal{E}]$. 
The total degree of a monomial $\prod_{E\in \mathcal{E}} E^{n_E}$, with $n_E$ a non-negative integer for all $E$, is the sum of the degrees of  the  variables, $\sum_{E\in \mathcal{E}} n_E$. The degree of a polynomial is   the maximum of the total degrees of its monomials. 

If a polynomial $p$ vanishes for all assignments  $a\colon\mathcal{E}\rightarrow \R_+$ then $p=0$ identically. 
Further, if $p$ is a non-zero polynomial  in $\R[\mathcal{E}]$ such that the degree of each variable in each monomial is either $1$ or zero, then all the coefficients of $p$ are  non-negative if and only if $p(a(\mathcal{E}))>0$ for any assignment  $a\colon\mathcal{E}\rightarrow \R_+$. If this is not the case then there is an assignment such that $p(a(\mathcal{E}))=0$.

For vectors $u=(u_1,\ldots,u_m)\in\R^m$ and $v=(v_1,\ldots,v_m)\in\R^m$, we let $u \wedge v$ be the component-wise minimum, $(u \wedge v)_i=\min(u_i,v_i)$, and let 
$$v_+ =(\max(v_1,0),\ldots,\max(v_m,0))\quad \textrm{and} \quad  v_-=(\min(v_1,0),\ldots,\min(v_m,0))$$ be the positive and negative parts of $v$. The support of $v$ is defined as the set of indices for which $v$ is non-zero, $\supp(v)=\{i\vert v_i\not=0\}$. The positive support of $v$ is $\supp^+(v)=\supp(v_+)$ and the negative support is $\supp^-(v)=\supp(v_-)$.  Let $v^t$ denote the transpose of $v$ and $u\cdot v$ the usual scalar product in $\R^n$.

For every $x\in \R$, we let $\sign(x)\in \{-,0,+\}$ be defined as
$$\sign(x)= \begin{cases} -  & \text{if }x<0, \\ 0 & \text{if }x=0, \\ +  & \text{if }x>0. \end{cases} $$
Signs are multiplied using the usual rules. If $\sigma$ is a sign and $x\in \R$ then $\sigma\cdot x$ is $0$ if $\sigma=0$ and $\pm x$ if $\sigma=\pm$, respectively.

We let $\# B$ denote the cardinality of a finite set $B$.

\section{Motivation: networks as dynamical systems}\label{crns}

In this section we introduce networks and  kinetics, and associate a dynamical system with a network and a kinetics.  The definition of a network is identical to that of a chemical reaction network,  which is used mainly in (bio)chemistry to describe networks of (bio)chemical reactions \cite{feinbergnotes}.
In general we use the nomenclature that is standard for chemical reaction networks. See for instance \cite{hornjackson,feinbergnotes,Feinbergss} for background and extended discussions.

\begin{definition}\label{crn}
A \emph{network} $\mmN$ consists of three finite sets:
\begin{enumerate}[(1)]
\item A set $\mmS=\{S_1,\dots,S_n\}$ of \emph{species}.
\item A set $\mmC\subset \omN^{n}$ of \emph{complexes}.
\item A set  $\mmR=\{r_1,\dots,r_m\}\subset \mmC\times \mmC$ of \emph{reactions}, such that $(y,y)\notin \mmR$ for all $y\in \mmC$, and if $y\in \mmC$, then there exists $y'\in \mmC$ such that either $(y,y')\in \mmR$ or $(y',y)\in \mmR$. 
\end{enumerate}
A network is denoted by $\mmN=(\mmS,\mmC,\mmR)$.
\end{definition}

We use the convention that  an element $r_j=(y_j,y'_j)\in \mmR$  is denoted  by $r_j\colon y_j\rightarrow y'_j$.
The \emph{reactant} and the \emph{product} (complexes) of a reaction $r_j\colon y_j\rightarrow y'_j$  are $y_j$ and $y'_j$, respectively. By definition,   any complex is either the reactant or the product  of some reaction.
The zero complex $0\in \mmC$ is allowed by definition. A reaction $S_i\rightarrow 0$, $S_i\in\mmS$, is called an \emph{outflow reaction}.

Throughout the paper, we use  $n$ to denote the number of species in  $\mmS$. 
The species $S_i$ is identified with the $i$-th canonical $n$-tuple of $\omN^n$ with $1$ in the $i$-th position and zeroes elsewhere. Accordingly, a complex $y\in \mmC$ is given as $y=\sum_{i=1}^n y_i S_i$ or  $(y_1,\dots,y_n)$. 
  We assume that $y\in\omN^n$ as reactions typically involve integer numbers of species. However, the results presented in this paper hold generally for $y\in\omR^n$. In examples we will often use other letters than $S_i$ for species to ease the presentation.
  Generally, we use $i$ to denote a species index and $j$ to denote a reaction index.

 \begin{example}\label{futile} 
Enzyme biology provides a rich source of examples. For instance, consider the network with set of biochemical species $\mmS=\{S_1,S_2,S_3,S_4\}$, set of complexes $\{S_{1} + S_2,S_1+S_3,S_2,S_3,S_4 \}$ and reactions 

\vspace{0.1cm}
 \centerline{\xymatrix{
S_{1} + S_2 \ar@<0.3ex>[r]  & S_4  \ar@<0.3ex>[l]  \ar[r] & S_{1} + S_3  &
S_3 \ar[r]  & S_2.
}}

\vspace{0.1cm}
\noindent  That is, the reactions are $r_1\colon  S_1+S_2\rightarrow S_4$, $r_2\colon  S_4\rightarrow S_1+S_2$, $r_3\colon S_4\rightarrow S_1+S_3$ and $r_4\colon S_3\rightarrow S_2$. This network assumes  the Michaelis-Menten enzyme mechanism in which a substrate $S_2$ is modified into a substrate $S_3$ through the  formation of an intermediate  $S_4$ \cite{enz-kinetics}. The reaction is catalyzed by an enzyme $S_1$. The modification can be reversed via a direct demodification reaction.
\end{example}

Reactions in a network are schematic representations of dynamical processes. Over time the concentrations or abundances of the species in the network change as a consequence of the reactions. In order to describe the dynamical properties of the network we introduce  a kinetics (Definition~\ref{kin}) and  the species  formation rate function (Definition~\ref{spec-func}). The kinetics provides the  reaction rate for given species concentrations and the species formation rate  function the instantaneous change in the concentrations when considering simultaneously the rate of all reactions.
 
\begin{definition}\label{kin} 
 A \emph{kinetics} for a network  $\mmN=(\mmS,\mmC,\mmR)$ 
 is an assignment to each reaction $r_j\in\mmR$ of a rate function $K_{j}\colon\Omega_K\rightarrow \omR_+$, where  $\Omega_K$ is a set such that  $\R^n_+\subseteq\Omega_K\subseteq\omR_+^n$, $c\wedge d\in\Omega_K$ whenever $c,d\in\Omega_K$, and 
 $$K_{j}(c)\geq 0\quad \text{ for all }\quad  c\in\Omega_K.$$ 
A kinetics for a network $\mmN$ is denoted by $K=(K_1,\dots,K_m)\colon \Omega_K\rightarrow \omR^m_+$.
If $K_{j}$ is differentiable for all $j=1,\dots,m$ and $c\in\R^n_+$ then $K$  is said to be a \emph{differentiable kinetics}.
\end{definition} 

\begin{example}
Let $\mmN=(\mmS,\mmC,\mmR)$  be the network with $\mmS=\{S_1,S_2,S_3\}$, $\mmC=\{S_1+S_2,S_3\}$ and $\mmR$ given by the reaction $r_1\colon S_1+S_2\rightarrow S_3$. The kinetics $K\in\mmK(\mmN)$ defined by $K_{1}(c)=k  c_1/((\beta+c_1)c_2^{\alpha})$, where $k,\alpha,\beta$ are positive constants has $\Omega_K=\omR\times \R_+\times\omR_+$.
\end{example}

\begin{example}\label{mass-action} 
Kinetics commonly used in chemistry and biology are the so-called mass-action kinetics. These were introduced by Guldberg and Waage in the 19th century based on the ideal assumption that the rate of a reaction is proportional to the product of the concentrations of the reactant species \cite{enz-kinetics}.
Specifically,  each reaction $r_j\colon y_j\rightarrow y_j'$ is assigned a positive constant $k_{j}\in \R_+$ and the rate function for the reaction is given by 
$$K_{j}(c)= k_{j} \prod_{i=1}^n c_i^{y_{j,i}}$$
with $\Omega_K=\omR^n_+$, $K=(K_1,\dots,K_m)$. We adopt the convention that $0^0=1$.
Under \emph{in vivo} conditions, however, the use of mass-action kinetics might not be fully justified. Reactant species might not form a homogeneous mixture,  for instance because they appear in too low concentrations or because their distribution depends on spatial constraints. 
In situations in which the use of mass-action is not justified, the use of other types of  kinetics  such  as \emph{power-law kinetics} or \emph{Hill-type kinetics} are often preferred. These will be introduced later. Note that $K_{j}(c)$ is an increasing function in $c_i$ and does not depend on other species than those involved in  $y_j$. 
\end{example}

The \emph{stoichiometric matrix} $A$ is defined as the $n\times m$ matrix whose $j$-th column is $y'_j-y_j$. 

\begin{definition}\label{spec-func}
The \emph{species  formation rate function} for a network $\mmN=(\mmS,\mmC,\mmR)$
with kinetics $K$ and stoichiometric matrix $A$ is the map $f_{A,K}\colon \Omega_K\rightarrow \Gamma$ defined by
$$f_{A,K}(c)=AK(c) = A\left(\begin{array}{c} K_{1}(c)  \\ \vdots \\ K_{m}(c) \end{array}\right)=\sum_{j=1}^m K_{j}(c)(y'_j-y_j).$$
\end{definition}

The dynamics of a network $\mmN$ with kinetics $K$ and stoichiometric matrix $A$ is described by a set of ordinary differential equations (ODEs) given by the species formation rate function: 
\begin{align} \label{ode}
\dot{c} &=f_{A,K}(c),
\end{align}
where $\dot{c}=\dot{c}(t)$ denotes the derivative of $c(t)$ with respect to $t$.
 Observe that  the image of $f_{A,K}$ is contained in $\im(A)$ and hence the dynamics of the system is confined to invariant linear spaces of the form $c+\im(A)$. 
In other words, for any $\omega\in \im(A)^{\perp}$ we have that $\omega\cdot \dot{c}=0$. Therefore, $\omega\cdot c$ is independent of time and determined by the initial concentrations of the system. The value of $\omega\cdot c$ is called a \emph{conserved amount}.

In this context, $\im(A)$ is called the \emph{stoichiometric space}. Two  vectors   $c,c'\in \R^n$   are called \emph{stoichiometrically compatible} if $\omega\cdot c = \omega\cdot c'$  for all $\omega\in \im(A)^{\perp}$, and $c,c'$ are said to be in the same \emph{stoichiometric class}, $c+\im(A)$. We let  $s$ be the rank of $A$ and, thus, the dimension of $\im(A)^{\perp}$ is $d=n-s$. 

The  \emph{steady states} of the network are the solutions to the system of  equations in $c_1,\dots,c_n$ obtained by setting the derivatives of the  concentrations to zero:
\begin{align*} 
0 =& f_{A,K}(c). 
\end{align*}
This system of equations is referred to as the steady-state equations. In particular, we are interested in the  positive steady states, that is, the  solutions $c$ to the steady-state equations such that all concentrations are positive, $c\in \R^n_+$.

\begin{example} 
The stoichiometric matrix of Example~\ref{futile} is
{\small 
\begin{equation}\label{stoich:futile}
A=\left( \begin{array}{rrrr}
          -1 & 1 & 1 & 0  \\ -1 & 1 & 0 & 1  \\ 0 & 0 & 1 & -1  \\ 1 & -1 & -1 & 0 
         \end{array}\right)
\end{equation}}
and has rank $s=2$.  A basis of  $\im(A)^{\perp}$  is
$\{ \omega^1,\omega^2\}$ with 
\begin{align}\label{perpbasis}
\omega^1 &= (1,0,0,1), & \omega^2 &= (0,1,1,1).
\end{align} 
If $K$ is any kinetics, then the  corresponding  system of ODEs is: 
{\small $$
\left(\begin{array}{c} \dot{c_1}\\  \dot{c_2} \\ \dot{c_3} \\  \dot{c_4} \end{array}\right) = \left( \begin{array}{rrrr}
          -1 & 1 & 1 & 0  \\ -1 & 1 & 0 & 1  \\ 0 & 0 & 1 & -1  \\ 1 & -1 & -1 & 0 
         \end{array}\right)\left(\begin{array}{c} K_{1}(c)  \\ K_{2}(c) \\ K_{3}(c) \\ K_{4}(c) \end{array}\right)  = \left(\begin{array}{c}  - K_{1}(c) +K_{2}(c) +K_{3}(c) \\ - K_{1}(c)+K_{2}(c)+K_{4}(c) \\  K_{3}(c) -K_{4}(c) \\  K_{1}(c) -K_{2}(c) -K_{3}(c) \end{array}\right).
$$}
Observe that $\dot{c}_1+\dot{c}_4=\dot{c}_2+\dot{c}_3+\dot{c}_4=0$ for any kinetics $K$. 
\end{example}

\section{Dynamical systems admitting a factorization}\label{sec:factor} 
Dynamical systems arising from reaction networks in the way specified above have a specific form, that is, the species formation rate function factors as the product of a matrix $A$ and a function vector $K\colon \Omega_K\rightarrow \omR^m_+$. 

Any dynamical system 
$ \dot{c} = f(c)$
such that $f(c)$ admits a factorization of the same form, $f(c)=AK(c)$, can be interpreted as arising from a network with stoichiometric matrix $A$. 
The reactions are however not uniquely given by $A$ and might not have a physical interpretation.

 \begin{example}\label{s-system} 
Savageau  \cite{savageau1988} considers a model of microbial growth
with ODE system 
\begin{align*}
\dot{c}_1 &= \alpha_1 c_1c_2 c_3^{-1} -  c_1, & \dot{c}_2 &=  c_4 - \beta_2 c_1c_2 c_3^{-1}, \\  
 \dot{c}_3 &=  c_4 - \beta_2 c_1c_2 c_3^{-1}, & \dot{c}_4 &=  \beta_2 c_1c_2 c_3^{-1}-  c_4,
\end{align*}
where  $\alpha_1,\beta_2>0$. It can be written as
{\small $$
\left(\begin{array}{c} \dot{c_1}\\  \dot{c_2} \\ \dot{c_3} \\  \dot{c_4} \end{array}\right) = \left(\begin{array}{cccc}  0 & 0 & 1 & -1 \\  -1 & 1 & 0 & 0 \\ -1 & 1 & 0 & 0 \\ 1 & -1 & 0 & 0    \end{array} \right) \left(\begin{array}{c}  \beta_2 c_1c_2 c_3^{-1}  \\ c_4 \\  \alpha_1 c_1c_2 c_3^{-1} \\ c_1 \end{array}\right).
$$}
It can, for example,  be interpreted as a network with reactions $S_2+S_3\rightleftharpoons S_4, 0\rightleftharpoons S_1$.  
 \end{example}
 
We study the steady states of dynamical systems of the form $\dot{c}=AK(c)$ for a fixed $n\times m$ matrix $A$ but with a variable kinetics $K(c)$. In order to simplify the notation, we use $f_K$ to denote $f_{A,K}$  whenever there is no confusion. 
The dynamics of the system $\dot{c}=AK(c)$ takes place in an invariant stoichiometric class given by the initial concentrations of the system. Hence, the relevant dynamical properties, such as multistationarity,  need to be inspected inside each stoichiometric  class. 

The main examples come from reaction networks and we keep the nomenclature introduced in the previous section. Even though the physical interpretation might be vague, we call $K$ a kinetics, $A$ the stoichiometric matrix and use the definitions introduced in the previous section.
We let $\mmK_{m,n}$ denote the set of kinetics $K\colon \Omega_K\rightarrow \omR^m_+$ for some $\Omega_K\subseteq \omR^n_+$ and let $\mmK_{m,n}^d\subset \mmK_{m,n}$ be the set kinetics that are continuous on $\Omega_K$ and differentiable on $\R^n_+$.

For any differentiable function $f=(f_1,\dots,f_q)\colon\Omega\rightarrow \R^q$ defined on a set $\Omega$ including $\R_+^m$, we let $J_c(f)$ denote the Jacobian of $f$ at $c\in\R^m_+$, that is, the $q\times m$ matrix  with entry $(j,i)$ being $\partial f_j(c)/\partial c_i$.  
If $K$  is a differentiable kinetics, then the Jacobian matrix $J_c(f_{K})$ factorizes as the product of two matrices:
 \begin{equation}\label{factor} J_c(f_{K}) = A (\partial K), \end{equation}
where $(\partial K)=(\partial K)(c)$ is the $m\times n$ matrix such that $(\partial K)_{j,i}= \partial K_{j}(c)/\partial c_i$.

Graphical conditions on the preclusion of multistationarity for dynamical systems for which the Jacobian can be decomposed as the product of two matrices (not necessarily of the form given here) have been studied in \cite{banaji-craciun1}.

\section{Some matrix-theoretical results}\label{matrixresults}
Let $M$ be an $n\times n$ matrix  and let $F$ be an $s$-dimensional  vector space that contains the space generated by the columns of $M$.  Let $F^{\perp}$ be the space orthogonal to $F$, which has dimension $d=n-s$.

\begin{definition}\label{extended}
A basis $\{\omega^1,\dots,\omega^d\}$  of  $F^{\perp}$ with $\omega^i=(\lambda_1^i,\dots,\lambda_n^i)$ is said to be \emph{reduced} if  $\lambda_i^i=1$ for all $i$ and $\lambda^i_j=0$ for all $j=1,\dots,\widehat{i},\dots,d$.
\end{definition}

After reordering of the coordinates of $\R^n$, if necessary, such a basis always exists and is unique.
Let $\widetilde{M}$ be the $n\times n$ matrix whose top $d$ rows are $\omega^1,\dots,\omega^d$ and the bottom $s$ rows agree with the bottom $s$ rows of $M$. We view $M$ as a linear map from $\R^n$ to $\R^n$ and let $ \ker(M)$ be the kernel of this map.

\begin{proposition}\label{kernel}
Let $M$, $F$ be as above. Let $\{\omega^1,\dots,\omega^d\}$ be a reduced basis of $F^{\perp}$ and $\wm$  the corresponding matrix. Then
$$ \ker(M)\cap F = \{0\}  \quad \text{if and only if}\quad \det(\widetilde{M})\neq 0.$$
\end{proposition}

For  any $n\times m$ matrix $B$, and sets $I\subseteq \{1,\dots,n\}$ and $J\subseteq\{1,\dots,m\}$, we let $B_{I,J}$ denote the submatrix of $B$ with entries of $B$ with indices $(i,j)$ in $(I,J)$.

\begin{proposition}\label{decomp} Let $M$ be an $n\times n$ matrix. Using the notation above, we have 
$$\det(\widetilde{M}) = \sum_{I\subseteq \{1,\dots,n\}, \#I=s} \det(M_{I,I}).$$
\end{proposition}
If  $M\subsetneq F$, then both sides of the equality are zero, because the rank of $M$ is strictly smaller than $s$.
In our applications, $F$ will be $\im(A)$ and $M$ will be $J_c(f_{K})=A(\partial K)$. In this case, using the Cauchy-Binet formula on the minors of a product of matrices, we have that 
\begin{equation}\label{cauchy}
\det(\widetilde{M}) = \sum_{I,J\subseteq \{1,\dots,n\}, \#I=\#J=s} \det(A_{I,J})\det( (\partial K)_{J,I}).
\end{equation}

\section{Degenerate steady states and injectivity}\label{sec:inj}

In this section some key concepts and definitions are introduced, namely that of degeneracy of a steady state  and  injectivity of a matrix.  
Denote the components of the species  formation rate function by $f_K=(f_{K,1},\ldots,f_{K,n})$.
Note that after reordering of the rows of a matrix $A$, if necessary, a reduced basis of $\im(A)^{\perp}$ always exists and is unique.  Therefore, from now on, we assume that the rows  of $A$ are ordered such that a  reduced basis exists. 

\begin{definition}\label{extended2}
Let $A$ be an $n\times m$ matrix of rank $s$ and $\{\omega^1,\dots,\omega^d\}$ the reduced basis of  $\im(A)^{\perp}$. If $K\in\mmK_{m,n}$ is a kinetics then the \emph{associated extended rate function} $\widetilde{f}_K\colon \Omega_K\rightarrow \R^n$ is the function defined by  
$$\widetilde{f}_K(c)=(\,\omega^1\cdot c,\dots,\omega^d\cdot c,f_{K,d+1}(c),\dots,f_{K,n}(c)\,).$$
\end{definition}

\begin{example} 
The basis of $\im(A)^{\perp}$ provided for Example~\ref{futile} in equation \eqref{perpbasis}  is reduced.
The associated extended rate function $\widetilde{f}_{K}\colon \Omega_K\rightarrow \R^4$ for any kinetics $K$ is
$$
\widetilde{f}_{K}(c) =(c_1+c_4,c_2+c_3+c_4,  K_{3}(c) -K_{4}(c), K_{1}(c) -K_{2}(c) -K_{3}(c)). $$
\end{example}

\begin{definition}\label{def:deg} Let $\dot{c}=AK(c)$ be a dynamical system such that  $K$  is a  differentiable  kinetics. A steady state $c\in \R_+^n$ of the system  if \emph{degenerate} if $\ker(J_{c}(f_{K}))\cap \im(A) \neq \{0\}$.
\end{definition}

That is, a steady state is degenerate if the Jacobian restricted to the stoichiometric subspace $\im(A)$ is non-singular. 
Using the constructions in Section~\ref{matrixresults} with $F=\im(A)^{\perp}$,
we have that 
\begin{equation}\label{jaccomp}
 \ker(J_{c}(f_{K}))\cap \Gamma = \{0\} 
\quad \text{if and only if}\quad \det(J_{c}(\widetilde{f}_{K}))\neq 0.
\end{equation}
It follows that a steady state $c\in \R^n_+$  is degenerate if and only if $\det(J_{c}(\widetilde{f}_{K}))=0$. 
The Jacobian of $\widetilde{f}_K$ has a natural interpretation as the flow of the dynamical system projected onto the stoichiometric space \cite{helton:determinant}.

Finally, we introduce the notion of injectivity. 
\begin{definition}\label{def-standard} 
Let   $A$ be an $n\times m$ matrix and $\mmK_0\subseteq \mmK_{m,n}$.
\begin{enumerate}[(i)]
 \item We say that  $A$ is \emph{injective} over $\mmK_0$ if for any pair of distinct stoichiometrically compatible vectors $a,b\in\R^n_+$ we have $AK(a)\not=AK(b)$ for all $K\in\mmK_0$.
\item The matrix $A$ is said to have the  \emph{capacity for multiple positive steady states} over $\mmK_0$ if there exists a kinetics $K\in\mmK_0$ and distinct  stoichiometrically compatible vectors  $a,b\in \R^n_+$ such that $AK(a)=AK(b)=0$. 
\end{enumerate}
\end{definition}

Note that being injective is equivalent to requiring that the function $\widetilde{f}_K$ is injective over  $\R^n_+$ for all  $K\in\mmK_0$.
Clearly, if $A$ is injective  over  $\mmK_0$, then $A$  does not have the capacity for multiple steady states over $\mmK_0$ and, thus, (i) implies (ii) in Definition \ref{def-standard}.

\begin{remark}
If $A$ is the stoichiometric matrix of a network $\mmN$, we say that the network $\mmN$ is injective and that the network $\mmN$ has the capacity for multiple steady states, if it is the case for $A$.
\end{remark}

The aim is to provide a criterion for a matrix $A$  to be injective  over a set of kinetics $\mmK_0$  in terms of computational tractable quantities.
To this end we introduce the class of power-law kinetics (defined in Section~\ref{gen-mass}) and derive some injectivity results for classes of power-law kinetics using  techniques introduced in \cite{Feliu-inj}.

\section{Power-law kinetics}
\label{gen-mass}

Power-law kinetics form a large family of kinetics \cite{hornjackson}.  They are generalizations of  mass-action kinetics and are based on a power-law formalism. Their  general form makes them flexible for modeling purposes  in many areas of science such as chemistry, ecology and epidemics. 
Furthermore, these kinetics are important in the  study of injectivity  in that they, in some sense, are   ``dense" in the set of kinetics that are strictly monotonic or constant in each concentration  (to be  made precise in Section~\ref{general-kinetics}). That is to say, injectivity of a matrix $A$ over certain sets of kinetics can be determined on the basis of injectivity of $A$ over suitable sets of power-law kinetics.

In this section we introduce power-law kinetics.   Dynamical systems with power-law kinetics have properties that are very similar to those with mass-action kinetics and  similar results regarding injectivity can be derived for power-law kinetics as for mass-action kinetics themselves.

For   a concentration vector $c$ and $v\in\R^n$, we associate the power law  $c^v=\prod_{i=1}^n  c_i^{v_i}$.   For example, if $v=(2.1,0.72,0,-1)\in \R^4$, then the associated power law is $c^v=c_1^{2.1}c_2^{0.72}c_4^{-1}$.  A power law is well defined for all $c\in\omR^n_+$ such that $c_i>0$ if $v_i<0$.
 
\begin{definition}\label{def-pl} 
A kinetics $K\in\mmK_{m,n}$ is a \emph{power-law kinetics} if 
$$ K_{j}(c)=k_{j}c^{v_{j}},\qquad \textrm{for } j=1,\dots,m,$$
with $k_j\in \R_+$ and $v_j\in \R^n$. Here $\Omega_K$ is the maximal common domain of definition of the power laws $c^{v_{j}}$, $j=1,\dots,m$, in the positive orthant.

Thus, a power-law kinetics is defined by an $m\times n$ matrix $V$ with rows $v_j$, $j=1,\ldots,m$. The matrix $V$ is called the \emph{kinetic order} and   $\kappa=(k_1,\dots,k_m)$  the \emph{rate vector}. The scalar $k_j$ is the \emph{rate constant} of reaction $j$.
For simplicity, we identify the pair  $(\kappa,V)$ with the kinetics $K$. 
We let $\mmK_{m,n}^g[V]$ denote the set of power-law kinetics $(\kappa,V)$ with arbitrary $\kappa$ but fixed $V$.
\end{definition}

Note that $k_{j}$ is a scalar while $v_{j}$ is a vector. 
 By definition,  power-law kinetics are differentiable kinetics.
 Given a network $\mmN=(\mmS,\mmC,\mmR)$ with $n$ species and $m$ reactions,  mass-action kinetics are  special types of power-law kinetics obtained  by considering the kinetic order with $v_{j}=y_j$ for all $j$.   Therefore, if we let $Y=(y_1,\dots,y_m)$ then the set of mass-action kinetics for a network $\mmN$ is $\mmK_{m,n}^g[Y]$.

\begin{example}\label{ex:det2}  
Consider Example \ref{futile}  with the reduced basis given in \eqref{perpbasis}. A kinetic orders is given as:
\begin{equation}\label{z:futile} 
 V=\left( \begin{array}{cccc}
          1 &  1 & 0 & 0 \\ 0 & 0 & 0  & 1   \\ 0 & 0 &  0 & 1   \\ v &  0 & 1 & 0          \end{array}\right),\qquad v\geq 0. 
\end{equation}  
When $v=0$, the kinetic order corresponds to mass-action.
If $v>0$, then the kinetic order $V$ accounts for the hypothetical fact that 
 the concentration of the modification enzyme $S_1$ acts as an enhancer or inhibitor of the demodification of $S_3$ to $S_2$, that is, of reaction $r_4$. 
In the latter case, the kinetics is:
\begin{align*}
 K_{1}(c) &= k_1c_1c_2,  & K_{2}(c)&=k_2 c_4, & K_{3}(c) &=  k_3c_4, & K_{4}(c) &= k_4c_1^vc_3.
\end{align*}
If for instance $v=0.5$, then $ K_{4}(c) = k_4c_1^{0.5}c_3$ while if $v=-0.5$, then $ K_{4}(c) = k_4c_1^{-0.5}c_3$.
\end{example}

\begin{example} 
After a suitable change of coordinates, ODE models of electrocatalytic oxidation of formic acid exhibit a power-law structure with a negative exponent. See for instance \cite[Examples 1,2]{sensse:analytic:2005}.
\end{example}

\begin{example}\label{lac}   
S-systems provide a rich source of examples of modeling with power-law kinetics. Consider the model of the lac gene circuit developed by Savageau  \cite{savageau:lac:2011}. The model has $5$ variables, $c_1,\dots,c_5$ and the ODEs take the form
\begin{align*}
\dot{c}_1 &= \alpha_1 c_4^{v_1} - \beta_1 c_1^{v_2}, & \dot{c}_2 &= \alpha_2 c_1^{v_3} - \beta_2 c_2^{v_4},& \dot{c}_3 &= \alpha_3 c_2^{v_5} - \alpha_4 c_2^{v_6}c_3^{v_7},\\  \dot{c}_4 &=  \alpha_4 c_2^{v_6}c_3^{v_7}- \alpha_5 c_2^{v_8}c_4^{v_9}, & \dot{c}_5 &=   \alpha_5 c_2^{v_8}c_4^{v_9}- \alpha_7 c_5^{v_{10}},
\end{align*}
for some positive exponents $v_*$ and positive constants $\alpha_*,\beta_*$ (equation (1) in \cite{savageau:lac:2011}). Similarly to Example~\ref{s-system}, this ODE system  factorizes as $AK$ (see also Example~\ref{lac2}).
\end{example}

\begin{example} 
So-called SIR (S=susceptible, I=infected, R=recovered) models are standard in epidemiology to describe the outbreak of an epidemics in a population.  One particular SIR model \cite{anderson-may} considers the network with set of species $\{S,I,R\}$ and reactions
$$r_1\colon S +I \rightarrow   2I,   \qquad r_2\colon I \rightarrow R. $$
The first reaction says that a susceptible individual might become infected in the presence of an infected. The second reaction says that infected individuals eventually recover. The SIR dynamics  can be expressed in  different ways. One possibility is the following set of differential equations \cite{stroud}:
$$
\dot{c}_1=-k_1 c_1^{\nu} c_2, \qquad 
 \dot{c}_2=k_1 c_1^{\nu} c_2 - k_2 c_2, \qquad
  \dot{c}_3= k_2 c_2,
$$ 
where $c_1,c_2,c_3$ are the concentrations of the species $S, I, R$, respectively, and $k_1,k_2>0$, $\nu> 0$ are the parameters of the model. That is, 
$K_{1}(c)= k_1 c_1^{\nu} c_2$, and $K_{2}(c)= k_2 c_2.$
 The parameter $\nu$ accounts for inhomogeneity in population mixing. If $\nu=1$ then the population is homogenous and the disease spreads according to the law of mass-action, whereas if $\nu\neq 1$ then the population is inhomogenous and the kinetics is a power-law kinetics.   If $\nu$ is allowed to be  negative, then susceptible individuals are repelled by infected individuals. Further, if we allow $\nu=0$, then the spread of the disease would be independent of the number of infected individuals. 
 \end{example}

\section{Injectivity for networks taken with power-law kinetics}\label{sec:powerlawinj}
In this section we provide criteria for a matrix to be injective with respect to the set of power-law kinetics with fixed kinetic order.
If  $K=(\kappa,V)\in\mmK_{m,n}$ is a power-law kinetics then the species  formation rate function is denoted by $f_{K}=f_{\kappa,V}$. If $\kappa$ is not fixed then the function $f_{\kappa,V}(c)$ can be seen as a polynomial function in the variables $k_{j}$ for all $j$.

\begin{theorem}\label{injecclose-pl} Let $A$ be an $n\times m$ matrix and $V$  an $m\times n$ kinetic order. Then the following are equivalent:
\begin{enumerate}[(i)]
\item $A$ is injective over $\mmK_{m,n}^g[V]$.
\item $\det(J_{c}(\widetilde{f}_{\kappa,V}))\neq 0$  for all $c\in \R^n_+$ and $\kappa\in\R_+^{m}$.
\end{enumerate}
\end{theorem}

The proof  is similar to the one given in \cite[Th. 5.6]{Feliu-inj} for mass-action kinetics and it is thus omitted here.
The following proposition provides an explicit description of $\det(J_{c}(\widetilde{f}_{\kappa,V}))$.  It is proven for mass-action kinetics in \cite{Feliu-inj} using a different approach. 
The current proof is based on the general matricial results explained in Section~\ref{matrixresults}.

\begin{proposition}\label{coefs2} Let $A$ be an $n\times m$ matrix of rank $s$ and  $V$  an $m\times n$ kinetic order. The determinant $\det(J_{c}(\widetilde{f}_{\kappa,V}))$ is a homogeneous polynomial in $\kappa=(k_1,\dots,k_m)$ of total degree $s$ and linear in each rate constant $k_{j}$. 

Further, let $J\subseteq \{1,\dots,n\}$ of cardinality $s$. The coefficient of the monomial $\prod_{j\in J} k_{j}$ in $\det(J_{c}(\widetilde{f}_{\kappa,V}))$ for $c\in \R^n_+$ is  
$$ c^{-\mathbf{1} + \sum_{j\in J} v_{j}} \sum_{I\in \{1,\dots,n\},\#I=s} \det(A_{I,J})\det(V_{J,I})  \prod_{i\notin I} c_{i}.$$
\end{proposition}

\begin{remark}
The determinant is in general \emph{not} a polynomial in  $V$ or in $c$ because the coordinates of  $V$ enter the expressions as exponents of $c$. However, the determinant $\det(V_{J,I})$ is a polynomial function in the non-zero  coordinates of $v_{j}$, $j\in J$,  excluding the entries with indices not in $I$.
\end{remark}

\begin{proposition}\label{det-criterion}
Let $A$ be an $n\times m$ matrix of rank $s$ and  $V$  an $m\times n$ kinetic order. The following are equivalent:
\begin{enumerate}[(i)]
\item $A$ is injective over $\mmK_{m,n}^g[V]$.
\item The non-zero products $\det(A_{I,J})\det(V_{J,I})$ have the same sign for all sets $I,J\subseteq\{1,\dots,n\}$ of cardinality $s$. Further $\det(A_{I,J})\det(V_{J,I})\neq  0$ for at least one choice of $I$ and $J$.
\end{enumerate}
\end{proposition}

\begin{example}\label{ex:det}
Consider Example \ref{futile} and the kinetic order introduced in Example \ref{ex:det2}.
We have 
$$
 \det(J_{c}(\widetilde{f}_{\kappa,V}))=(k_2+k_3)k_4c_1^vc_4+k_1k_3c_1c_4+k_1k_4c_1^v(c_1+vc_3+c_2).
$$
 If $v\geq 0$, then all the terms of the determinant expansion have the same sign and by Theorem~\ref{injecclose-pl}  the network is injective over $\mmK_{m,n}^g[V]$. If, on the contrary, $v<0$, then the term $t=vk_1k_4c_1^vc_3$ has sign opposite to the rest of the terms. It follows that the network is not injective over $\mmK_{m,n}^g[V]$ if $v<0$. 
  The term  $t$ corresponds to the sets $J=\{1,4\}$ and $I=\{1,2\}$. Indeed, for these  sets we have, see \eqref{stoich:futile} and \eqref{z:futile}, 
$$ V_{J,I} = \left(\begin{array}{cc} 1 & 1  \\  v & 0 \end{array}\right),\qquad 
A_{I,J} = \left(\begin{array}{cc} -1 & 0 \\  -1 & 1  \end{array}\right),  $$
so that $ \det(V_{J,I})=-v$, $\det(A_{I,J})=-1$ and 
$c^{-\mathbf{1} + \sum_{j\in J} v_j}\prod_{i\notin I} c_{i} =(c_1^vc_4^{-1})(c_3c_4)=c_1^vc_3$. The sign of $t$ depends on that of $v$, unless  $v=0$, in which case the term vanishes.
\end{example}

\begin{example}\label{lac2}
Consider Example~\ref{lac}. 
The matrices $A$ and $V$ are
{\small $$ A=\left(\begin{array}{cccccccc}  1 & -1 & 0 & 0 & 0 & 0 & 0 & 0 \\
0& 0 & 1 & -1 & 0 & 0 & 0 & 0 \\
0 & 0 & 0 & 0 & 1 & -1 & 0 & 0 \\
0 & 0 & 0 & 0  & 0 & 1 & -1 & 0   \\
0 & 0 & 0 & 0 & 0 & 0 & 1 & -1
 \end{array}\right),\quad 
 V=\left(\begin{array}{ccccc}  
 0 & 0 & 0 & v_1 & 0 \\  v_2 & 0 & 0 & 0 & 0 \\  v_3 & 0 & 0 & 0 & 0 \\  0 & v_4 & 0 & 0 & 0 \\
  0 & v_5 & 0 & 0 & 0 \\  0 & v_6 & v_7 & 0 & 0 \\  0 & v_8 & 0 & v_9 & 0 \\ 0 & 0 & 0 & 0 & v_{10}
 \end{array}\right).$$}
The matrix $A$ has maximal rank. The products $\det(A_{I,J})\det(V_{J,I})$ are non-zero for the pairs $(I,J)$ with $I=\{1,2,3,4,5\}$ and 
$J=\{1,3,5,6,8\}, \{1,3,6,7,8\}$ or  $\{2,4,6,7,8\}.$ For these pairs, the determinant products are
$$ -v_{1}v_{3}v_{7}v_{8}v_{10},\qquad  v_{1}v_{3}v_5v_{7}v_{10},\qquad -v_{2}v_{4}v_{7}v_{9}v_{10},$$
respectively. Since the terms do not have the same sign, it follows that for any choice of positive exponents $v_{\ell}$ the network is not injective.
\end{example}

\section{Influence specifications}\label{sec:int}
In the previous section we  studied injectivity of the system $\dot{c}=AK(c)$, where $K$ is a power-law kinetics with fixed kinetic order and varying rate constants. In the following sections we will study injectivity  when $K$ is not fixed but belongs to a general family of kinetics. The family is given by requiring that $K$ fulfills certain monotonicity constraints, which are encoded by the sign-pattern of a matrix. When $K$ is differentiable, the sign-pattern agrees with the sign-pattern of the Jacobian of $K$.

\begin{definition}
\label{infl-spec}
An $m\times n$ \emph{influence specification}  is a sign matrix $Z=(z_{j,i})$, that is, a matrix whose entries are the signs $+,-,0$.  
\end{definition}

Let  $Z$ be an $m\times n$ influence specification.  For $j=1,\dots,m$, define
\begin{align*}
z_j^+&=\{i|\ z_{j,i}=+\}, &  z_j^-& =\{i|\ z_{j,i}=-\}, &  z_j^0&=\{i|z_{j,i}=0\}.
\end{align*}
Two concentration vectors  $a,b\in \omR_+^n$ are said to be \emph{non-overlapping with respect to  $Z$}, 
  $$\text{if }\quad  z_j^+\nsubseteq \supp(a)  \quad \text{ implies }\quad  z_j^+\subseteq \supp(b)\quad \textrm{for all }j. $$
That is, the coordinates $a_i,b_{\ell}$ (potentially with  $i=\ell$) cannot both be zero if $i,\ell$ are both in $z_j^+$ for some $j$.   When it is clear from the context what influence specification we are referring to, we omit ``with respect to $Z$" and just say that $a,b$ are non-overlapping.
Non-overlapping is a concept that specifies how two concentration vectors on the  boundary of the positive orthant are positioned with respect to each other in relation to an influence specification $Z$. In particular, if one or both of the concentration vectors $a,b$   are positive then they are non-overlapping. 

\subsection{Strictly monotone kinetics} 

We start with a definition.
\begin{definition}\label{respect}
A  kinetics $K\in \mmK_{m,n}$  \emph{respects an $m\times n$ influence specification $Z$} if, for all $c\in\Omega_K$,
$$K_{j}(c)>0\quad \text{ if and only if } \quad z_j^+\subseteq \supp(c)\quad \textrm{for all }j.$$
\end{definition}

Let $\Omega_K(z_j^+):=\{c\in \Omega_K|\ z_j^+\subseteq \supp(c)\}=\{c\in \Omega_K|\ K_{j}(c)\neq 0\}$ denote the set of concentration vectors for which the kinetics $K_{j}$ does not vanish.

\begin{definition}\label{def-mono0}
Let  $K\in \mmK_{m,n}$ be a kinetics  that respects an $m\times n$ influence specification $Z$. We say that $K$  is \emph{strictly monotonic with respect to} $Z$ if for all $j=1,\dots,m$ and  $i=1,\ldots,n$,
the restriction of the function  $K_{j}(\cdot)$ to $\Omega_K(z_j^+)$ is
\begin{enumerate}[(i)]
\item strictly increasing in $c_i$ if $i\in z_j^+$.
\item strictly decreasing in $c_i$ if $i\in z_j^-$.
\item constant in $c_i$  if $i\in z_j^0$.
\end{enumerate}
Let $\mmK_{m,n}(Z)$ denote the set of   kinetics that are strictly monotonic with respect to the  influence specification $Z$. 
\end{definition}

The definition says that the rate functions $K_{j}$ are strictly monotonic or constant in the coordinate $c_i$ whenever the remaining coordinates take positive values for the species with positive influence. 
Strictly monotonic refers to the kinetics, but this does not imply that the species formation rate function is monotone.

\begin{example} \label{larvae}
The length of a larvae is often  assumed to increase linearly with a slow down in the growth rate as the length increases. 
Denoting by $c$ the length of the larvae, one model of the length is $\dot{c}=\alpha_1(\alpha_2+c),$
where $\alpha_1,\alpha_2$ are positive constants. 
In this system, $n=m=1$ and $A=(1)$. The kinetics  $K_{1}(c)=\alpha_1/(\alpha_2+c)$ 
 is strictly monotonic with respect to the influence specification given by 
 $Z=(-)$.
\end{example}

The following lemma shows how Definition~\ref{def-mono0} can be stated in  the terminology of \cite{shinar-conc}. We will use this characterization again in Section~\ref{sec:shinar}.

\begin{lemma}\label{lemma-mono0}
Let  $K\in \mmK_{m,n}$ be a kinetics  that respects an $m\times n$ influence specification $Z$. Then, $K$ is strictly monotonic with respect to  $Z$ if and only if  for each pair of non-overlapping concentration vectors $a,b\in\Omega_K$, the following implications hold for all $j$: 
\begin{enumerate}[(i)]
\item if $K_{j}(a)>K_{j}(b)$ then $\sign(a_i-b_i)=z_{j,i}\neq 0$ for some $i$.
\item if $K_{j}(a)=K_{j}(b)$ then either $a_i=b_i$ for all  $i\in z_j^+\cup z_j^-$, or  $\sign(a_i-b_i)=z_{j,i}\neq 0$ and $\sign(a_{\ell}-b_{\ell})=-z_{j,\ell}\not=0$ for some distinct $i,\ell$.
\end{enumerate}
\end{lemma}

\begin{example}\label{pop-growth}
A standard model of population growth is 
\begin{equation}\label{pop}
\dot{c}=r c\left(1-\frac{c}{D}\right)
\end{equation}
\cite{may}, where $c$ denotes the size of a population  $S$ and $r, D$ are  positive constants. 
The system factorizes as $AK$ with $A=(1, -1)$ and $K=(rc,rc^2/D)$. 
 The kinetics  is strictly monotonic with respect to $Z=(+,+)^t$. The system has a stable steady state at $c=K$ and an unstable steady state at $c=0$. The two steady states are non-overlapping with respect to $Z$.
\end{example}

\begin{definition}\label{def-diff}
A kinetics $K\in \mmK_{m,n}$  is \emph{differentiable with respect to an $m\times n$ influence specification} $Z$ if $K$ respects $Z$ and, for every $j$, $K_{j}(\cdot)$ is continuous at $c\in\Omega_K$,  differentiable at $c\in\R^n_+$, and for each index $i$ the partial derivative $\frac{\partial K_{j}}{\partial c_i}(c)$ has constant sign
$$z_{j,i}= \sign\left(\frac{\partial K_{j}}{\partial c_i}(c)\right)$$
 in $\R^n_+$.
Let $\mmK^{d}_{m,n}(Z)$ denote the set of   kinetics that are  differentiable with respect to the influence specification $Z$. 
\end{definition}

Note that $\mmK^{d}_{m,n}(Z)$ is \emph{not} the intersection of $\mmK_{m,n}(Z)$ with $\mmK_{m,n}^d$ as we require the sign of the partial derivatives to be constant in Definition~\ref{def-diff}, which is not implied by being strictly monotonic and differentiable. 
 We state without proof:

\begin{lemma}\label{diff-to-weak}
If $K\in\mmK^{d}_{m,n}(Z)$ then $K\in\mmK_{m,n}(Z)$. 
\end{lemma}

 \begin{example}\label{s-system2} Consider Example~\ref{s-system}. The kinetics is strictly monotone with respect to the influence specification given by
{\small $$Z=\left( \begin{array}{cccc} + & + & - & 0 \\ 0 & 0 & 0 & + \\  + & + & - & 0 \\ + & 0 & 0 & 0 \end{array} \right). $$}
The matrix $Z$ is simply the sign pattern of the Jacobian of $K$.
 \end{example}

\begin{example}\label{jacob-monod}
Jacob and Monod \cite{jacob-monod} 
proposed a model for  bacterial uptake of nutrients in microbial ecology. The modeling equations are 
$$\dot{c}_1= \frac{\alpha c_1 c_2}{\beta+c_2},\quad \text{ and } \quad \dot{c}_2=-\frac{ \alpha c_1c_2}{\gamma(\beta+c_2)},$$
where $\alpha,\beta,\gamma$
are positive parameters and $c_1,c_2$ are the concentrations of two species $S_1, S_2$. The species $S_1$ is a bacteria that feeds on a (chemical) nutrient $S_2$. The species formation rate function factorizes as $AK(c)$ with
$$A=\left(\begin{array}{cc} 0 & -1 \\1 & 0 \end{array}\right)\quad \textrm{ and }\quad   
K=\left(\frac{\alpha c_1 c_2}{\beta+c_2},\frac{ \alpha c_1c_2}{\gamma(\beta+c_2)}\right).$$ 
The kinetics is strictly monotonic with respect to the influence specification $Z$ with $+$ in all entries.
All steady states have $c_1=0$ or $c_2=0$. Hence two steady states are non-overlapping with respect to $Z$ if at least one of them is strictly positive. Further, all steady states are degenerate.
\end{example}

\subsection{Networks and influence specifications}
The term ``influence specification'' comes from  interpreting a dynamical system as a network $\mmN=(\mmS,\mmC,\mmR)$ with species $S_i$ and reactions $r_j$. 
Since the influence specification indicates the behavior of $K_j$ with respect to each variable,
the species $S_i$ with  $z_{j,i}\neq 0$ influence the reaction $r_j$. The species with $z_{j,i}=+$ are assumed to have  positive influence or enhance the  reaction, while those with $z_{j,i}=-$ are assumed to have negative influence and   an inhibitory effect on the reaction.  Those with $z_{j,i}=0$ have  neutral influence.
A reaction can only occur  if all species with positive influence on the reaction are present (that is, they are in positive concentrations). If one species is not present  then the reaction cannot occur.  Absence of species with negative or neutral influence does not prevent the reaction from taking place.

Given a network, we define the following distinguished influence specifications: 
\begin{itemize}
\item The \emph{complex dependent  influence specification}, denoted by $Z_{\mmC}$,  is defined as the influence specification with $z_{j,i}=+$ for  $i\in \supp(y_j)$ and zero otherwise.
\item The \emph{reaction dependent influence specification}, denoted by $Z_{\mmR}$,  is defined as the influence specification with $z_{j,i}=\sign(y_{j,i}-y'_{j,i})$.
\end{itemize}
 Note that  any mass-action kinetics $K$ belongs to $\mmK^{d}_{m,n}(Z_{\mmC})$.  

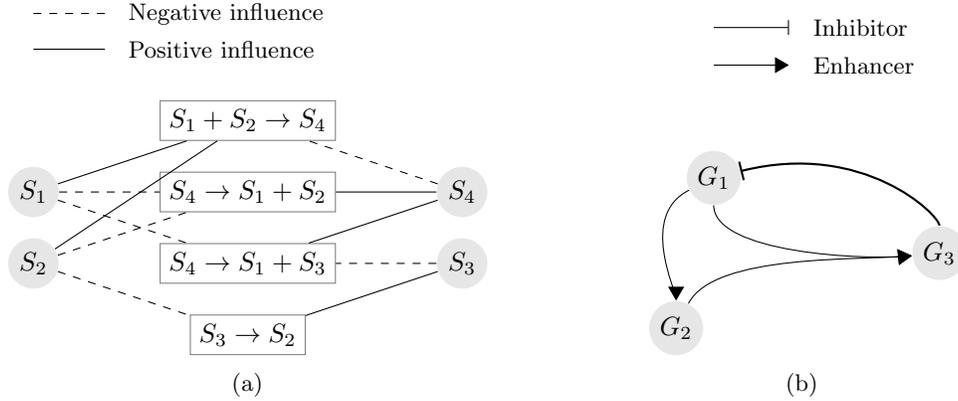
\begin{figure}[!t]
\centering
\subfigure[]{
{\small \begin{tikzpicture}[scale=0.95]
\node[rectangle,draw=gray] (r1) at (0,0) {$S_1+S_2\rightarrow S_4$};
\node[rectangle,draw=gray] (r2) at (0,-1) {$S_4\rightarrow S_1+S_2$};
\node[rectangle,draw=gray] (r3) at (0,-2) {$S_4\rightarrow S_1+S_3$};
\node[rectangle,draw=gray] (r4) at (0,-3) {$S_3\rightarrow S_2$};

\node[circle,fill=gray!20!white,inner sep=2pt] (s1) at (-3,-1) {$S_1$};
\node[circle,fill=gray!20!white,inner sep=2pt] (s2) at (-3,-2) {$S_2$};
\node[circle,fill=gray!20!white,inner sep=2pt] (s3) at (3,-2) {$S_3$};
\node[circle,fill=gray!20!white,inner sep=2pt] (s4) at (3,-1) {$S_4$};

\draw[-] (s1) -- (r1);
\draw[-,dashed] (s1) -- (r2);
\draw[-,dashed] (s1) -- (r3);
\draw[-,dashed] (s4) -- (r1);
\draw[-] (s4) -- (r2);
\draw[-] (s4) -- (r3);
\draw[-] (s3) -- (r4);
\draw[-] (s2) -- (r1);
\draw[-,dashed] (s2) -- (r2);
\draw[-,dashed] (s3) -- (r3);
\draw[-,dashed] (s2) -- (r4);

\draw[-,dashed] (-3,1.5) -- (-2,1.5);
\node[anchor=west] (t1) at (-1.8,1.5) {\footnotesize Negative influence};
\draw[-] (-3,1) -- (-2,1);
\node[anchor=west] (t1) at (-1.8,1) {\footnotesize Positive influence};
\end{tikzpicture} } \label{intgraph}}
\hspace{1.5cm}
\subfigure[]{{\small \begin{tikzpicture}
\node[circle,fill=gray!20!white,inner sep=2pt] (G1) at (-1,1) {$G_1$};
\node[circle,fill=gray!20!white,inner sep=2pt] (G2) at (-1.5,-1) {$G_2$};
\node[circle,fill=gray!20!white,inner sep=2pt] (G3) at (2,0) {$G_3$};

\begin{scope}[>=triangle 60] 
\draw[->] (G1) .. controls +(-1,-0.5) and +(0,0.5) .. (G2);
\draw[-|,thick] (G3) .. controls +(0,0.5) and +(2,0.5) .. (G1);
\draw[->] (G1) .. controls +(0,-1) and +(-1,-0.1) ..  (G3);
\draw[->] (G2) .. controls +(0.5,1) and +(-1,-0.1) ..  (G3);
\draw[-|] (-1,3) -- (0,3);
\node[anchor=west] (t1) at (0.2,3) {\footnotesize Inhibitor};
\draw[->] (-1,2.5) -- (0,2.5);
\node[anchor=west] (t1) at (0.2,2.5) {\footnotesize Enhancer};
\end{scope}
\end{tikzpicture} }\label{Xgenes}}

\caption{Graphical representation of a network with an influence specification. (a) Reaction dependent influence specification for Example \ref{futile}, drawn as a bipartite graph. There is an edge between a species and a reaction if the species has non-zero influence on the reaction. If the influence is positive, the edge is solid. If the influence is negative, we draw a dashed edge. The complex dependent influence specification is obtained by removing the dashed edges. The species interacting in a reaction, for example $S_4\rightarrow S_1+S_3$ cannot be read off from the edges. (b) The influence specification for the transcription of the three genes in Example \ref{gene-regu} are shown. The presence of $G_3$ reduces the production of $G_1$, whereas the presence of $G_1$ and $G_2$ cooperatively induce the production of $G_3$. Likewise $G_1$ induces the production of $G_2$.}
\end{figure}

It is sometimes useful to illustrate an influence specification with a labeled bipartite graph with node set 
$\mmS\cup \mmR$. We draw a positive edge  between a species and a reaction if the species has positive influence over the reaction.  We draw a negative edge if the species has negative influence on the reaction.  
Example~\ref{futile} is illustrated in Figure~\ref{intgraph}, assuming a reaction dependent influence specification.

\begin{example}\label{gene-regu} 
Common examples in the literature are gene regulatory networks \cite{karlebach}.
These are typically represented by diagrams as the one shown in Figure \ref{Xgenes}. The diagram represents 
three genes, each transcribing a protein. The proteins mutually affect the transcription rates of the genes such that the genes influence each other. The diagram corresponds to the  network  with production reactions
  $$ \xymatrix{0 \ar[r]  & G_i}\qquad i=1,2,3, $$
  and influence specification  $Z$  with $z_{1,3} = -$, $z_{2,1} = z_{3,1} = z_{3,2} = +,$ and zero otherwise. 
\end{example}

\begin{example}\label{lotka}  
An example from ecology is given by the Lotka-Volterra equations for modeling predator-prey dynamics or competing populations \cite{may,murray}. The modeling equations  are 
$$\dot{c}_1= c_1(\alpha-\beta c_2),\quad \text{ and } \quad \dot{c}_2=-c_2(\gamma-\beta c_1),$$
where $c_1$ and $c_2$ are the abundance of two species $S_1$ (prey) and $S_2$ (predator). The system comes from a network with reaction set $\mmR=\{S_1 \rightarrow 2S_1,S_1 + S_2\rightarrow 2S_2,S_2\rightarrow 0\}$ and kinetics
$$K_{1}(c_1,c_2)= \alpha c_1, \quad K_{2}(c_1,c_2)= \beta c_1c_2, \quad \text{ and } \quad K_{3}(c_1,c_2)=\gamma c_2,$$ 
with $\alpha,\beta,\gamma>0$. 
The kinetics is mass-action and thus belongs to $\mmK_{3,2}(Z_{\mmC})$.
The model has been proposed independently in epidemics \cite{kermack,anderson-may}, as well as a in physical chemistry as a model of H$_2$O$_2$ combustion \cite{semenov}.
\end{example}

\subsection{Kinetic orders and influence specification}
A kinetic order for a power-law kinetics is intimately related to an influence specification. 
If $V$ is a kinetic order, define an influence specification $Z(V)$ by:  
$$Z(V)_{j,i} =\sign(v_{j,i}).$$ 
In this case we say that $Z(V)$ is the influence specification associated to the kinetic order $V$.
Reciprocally, if $Z$ is an influence specification, let  $V(Z)$  be the kinetic order   defined by $$V(Z)_{j,i} = z_{j,i}\cdot 1.$$
Note that   $Z(V(Z))=Z$.

\begin{example}\label{ex:kineticorder} 
Consider Example \ref{futile} and the kinetic order $V$ introduced in Example \ref{ex:det2}. Then we have:
{\small $$
Z(V)= \left( \begin{array}{cccc}
      + &  + & 0 & 0 \\ 0 & 0 & 0  & +   \\ 0 & 0 &  0 & +   \\ \sign(v) &  0 & + & 0   \end{array}\right)\quad \textrm{and}\quad V(Z(V))= \left( \begin{array}{cccc}
      1 &  1 & 0 & 0 \\ 0 & 0 & 0  & 1   \\ 0 & 0 &  0 & 1   \\ \sign(v)\cdot 1 &  0 & 1 & 0   \end{array}\right).
$$}
\end{example}

For an influence specification $Z$, let  $\mmK_{m,n}^g(Z)\subseteq \mmK_{m,n}(Z)$ be the set of power-law kinetics that are strictly monotonic with respect to $Z$.  If $V$ is a kinetic order  then it is straightforward to see that any power-law kinetics $(\kappa,V)$ is differentiable with respect to the associated influence specification $Z(V)$.
Therefore, 
$\mmK^{g}_{m,n}(Z)\subseteq \mmK^d_{m,n}(Z)$.
 Likewise any kinetics $(\kappa,V)$ that is strictly monotonic with respect to an influence specification $Z$ fulfills $Z=Z(V)$. Hence, we have the following lemma:
\begin{lemma} \label{subsets}
For any kinetic order $V$ we have $\mmK_{m,n}^g[V]\subseteq \mmK^{g}_{m,n}(Z(V))$. A kinetics $(\kappa,V)$ belongs to $\mmK^{g}_{m,n}(Z)$ if and only if $Z(V)=Z$.  Further
$$ \mmK^{g}_{m,n}(Z) = \bigcup_{V| Z(V)=Z} \mmK^{g}_{m,n}[V].$$ 
\end{lemma}

\subsection{$Z$-injectivity}
In this subsection we introduce the notion of $Z$-injectivity.

\begin{definition}\label{inject-nonover}
Let $A$ be an $n\times m$ matrix, $Z$ an $m\times n$ influence specification and let $\mmK_0\subseteq \mmK_{m,n}(Z)$. We say that $A$ is  \emph{$Z$-injective} over $\mmK_0$ if for any $K\in\mmK_0$ and any pair of   stoichiometrically compatible vectors $a,b\in\Omega_K$, that are non-overlapping with respect to $Z$, we have $f_{K}(a)\not=f_{K}(b)$.
\end{definition}

If  $A$ is $Z$-injective, then the existence of pairs of distinct non-overlapping stoichiometrically compatible steady states $a,b\in\Omega_K$ is precluded. In particular,    some types of multiple steady states at the boundary of $\Omega_K$ are excluded. 
 Since pairs of positive concentration vectors are non-overlapping for any influence specification $Z$,  $Z$-injectivity implies injectivity. 
 Shinar and Feinberg \cite{shinar-conc} exclude the occurrence of pairs of distinct stoichiometrically compatible steady states $a,b$ such that at least one of them is in the interior of $\Omega_K=\omR^n_+$ (in their setting).  As noticed above, such pairs are non-overlapping and hence covered by our approach.  In  \cite{Feliu-inj}, the  condition  of non-overlapping is applied to networks with mass-action kinetics and influence specification $Z_{\mmC}$. It is possible to have $f_K(\ta)=f_K(\tb)$ for a pair of \emph{not} non-overlapping vectors and at the same time $f_K(a)\not=f_K(b)$ for all non-overlapping pairs \cite{Feliu-inj}.

The rest of the paper is devoted to characterize the matrices  that are $Z$-injective for different families of kinetics and that consequently cannot have the capacity for multiple positive steady states. The aim is to provide a determinant criterion for an $n\times m$ matrix  to be $Z$-injective  over a set of kinetics $\mmK_0\subseteq \mmK_{m,n}(Z)$  in terms of computational tractable quantities (Theorem~\ref{equiv}).
To this end we relate injectivity for power-law kinetics  to $Z$-injectivity for    kinetics that are strictly monotonic with  respect to $Z$.

\subsection{Kinetic orders with common influence specification}

In Section \ref{sec:powerlawinj} we  studied  injectivity of matrices for a fixed kinetic order $V$ and arbitrary rate vector $\kappa$. Here we are concerned about  injectivity of matrices over a set of power-law kinetics for which the kinetic orders are related through their associated influence specifications.

To proceed we need  some additional notation and definitions. 
We first introduce a partial order on the set of $m\times n$ influence specifications. One influence specification $\widetilde{Z}$ is said to be smaller than another influence specification $Z$ if $\widetilde{z}^+_{j} \subseteq z^+_{j}$ and $\widetilde{z}^-_{j} \subseteq z^-_{j}$ for all $j=1,\dots,m$. If this is the case then we write $\widetilde{Z}\preceq Z$. 
The minimal element in this order is the zero influence specification, that is, the influence specification is zero for all reactions on all species.
There is not a unique maximal element in this order and all maximal elements must  fulfill $z^+_{j}\cup z^-_{j}=\{1,\ldots,n\}$.

\begin{definition}\label{signed} 
Let $A$ be an $n\times m$ matrix of rank $s$, $Z$ be an $m\times n$ influence specification and let
$\Sigma(Z) := \{ V |\ Z(V)=Z\}. $
We say that $Z$ has a \emph{signed $A$-determinant}  if the function
$$\delta  \colon \Sigma(Z)  \rightarrow \{-,0,+\},\qquad V \mapsto \sign(\det(\widetilde{AV}))$$ is constant. If, further,  the image of  $\delta$ is not zero then $Z$ is called \emph{$A$-sign-nonsingular} (A-SNS).
\end{definition}

If $A$ is fixed,  we simply say that $Z$ has a signed determinant.
According to the definition, $Z$ has a signed determinant if the signs of $\det(\widetilde{AV})$ and $\det(\widetilde{AV'})$ agree for any two kinetic orders $V,V'$ with influence specification $Z$. 
The  kinetic orders  in $\Sigma(Z)$ have the same sign-pattern, hence a kinetic order $V$ in $\Sigma(Z)$ is uniquely identified by the absolute values $|v_{j,i}|$  of the nonzero entries of $V$. Let $|V|$ denote the matrix obtained from $V$ by considering the absolute value component-wise. Then $V=Z\ast |V|$ (where $\ast$ denotes component-wise sign-number product). It follows that $\Sigma(Z)$ admits a positive parameterization. 

Let $X=\{x_{*,*}\}$ be a generic symbolic $m\times n$ matrix and $Z_X:=Z\ast X$. Let $p_Z(X):=\det(\widetilde{AZ_X})$. Then $p_Z(X)$ is a polynomial in the entries of $X$ for which $z_{j,i}\neq 0$, such that
\begin{equation}\label{detpz}
\det(\widetilde{AV}) = p_{Z}(|V|).
\end{equation}
 In fact,  $p_{Z}(X)$ is either the zero polynomial or a homogeneous polynomial of degree $s$ in $x_{j,i}$. Further, the degree of each monomial in each variable $x_{j,i}$ is either zero or one.

If $Z'$ satisfies $Z'\preceq Z$, then $p_{Z'}(X)$ is obtained from $p_Z(X)$ be setting some variables $x_{j,i}$ to zero.

\begin{lemma}\label{support} Let $A$ be an $n\times m$ matrix.
\begin{enumerate}[(i)]
\item $Z$ has a signed A-determinant if and only if $p_{Z}$ is either the zero polynomial or the non-zero coefficients of the monomials of  $p_{Z}$ have common signs.
\item Let $Z'$ be   such that $Z'\preceq Z$. 
\begin{enumerate}[(a)]
\item If $Z$ has a signed A-determinant then $Z'$ has a signed A-determinant.
\item  If  $Z$ and $Z'$ are both A-SNS  then $\delta(Z)=\delta(Z')$.
\item If $Z'$ is A-SNS and $Z$ has a signed A-determinant then $Z$ is A-SNS.
\end{enumerate}
\end{enumerate}
\end{lemma}
 
\medskip
Using equation \eqref{cauchy}, we have the decomposition
\begin{equation}\label{cauchy2}
\det(\widetilde{AZ_X}) = \sum_{I,J\subseteq \{1,\dots,n\}, \#I=\#J=s} \det(A_{I,J})\det( (Z_X)_{J,I}). 
\end{equation}
The product $\det(A_{I,J})\det( (Z_X)_{J,I})$ is  a polynomial, $p_{Z,I,J}$, in  the entries of $X$, such that
\begin{equation}\label{pJC}
p_{Z}(X) = \sum_{I,J\subseteq \{1,\dots,n\}, \#I=\#J=s}p_{Z,I,J}(X).
\end{equation}
Each polynomial $p_{Z,I,J}$  involves different variables and, hence,  none of the terms cancel out in the sum over $I,J$ (unless they are zero).
We say that  $p_{Z,I,J}$ is \emph{sign-nonzero} if it is non-zero and the coefficients of the non-zero terms in the polynomial have the same sign. The sign of any of the coefficients is then the sign of  $p_{Z,I,J}$. Consequently, we have:

\begin{corollary}\label{polyterms}
Let  $A$ be an $n\times m$ matrix and $Z$ be an $m\times n$  influence specification. Then $Z$ is A-SNS  if and only if  the non-zero terms $p_{Z,I,J}$ in \eqref{pJC} are sign-nonzero, have the same sign, and at least one of the terms is non-zero. 
\end{corollary}

An $n\times m$ matrix $A$ is injective over $\mmK^g_{m,n}(Z)$ if and only if $A$ is injective over $\mmK^g_{m,n}[V]$ for all kinetic orders with $Z(V)=Z$ (Lemma~\ref{subsets}). Hence, we can use Proposition~\ref{det-criterion} and Lemma~\ref{support} to derive a determinant criterion valid for $\mmK^g_{m,n}(Z)$.   Similarly, if $Z_1,Z_2$ are influence specifications  such that $Z_1\preceq Z_2$  then we can  use  Proposition~\ref{det-criterion} and Lemma~\ref{support} to derive a determinant criterion for $A$ to be injective over $\bigcup_{Z \vert Z_1\preceq Z\preceq Z_2}\mmK^{g}_{m,n}(Z)$.

\begin{proposition}\label{openpos2}   
Let  $A$ be an $n\times m$ matrix and $Z$ an $m\times n$  influence specification.  Then  $A$   is  injective over  $\mmK^g_{m,n}(Z)$ if and only if $Z$ is A-SNS.
\end{proposition}

Using equation \eqref{cauchy}, we can rephrase the previous proposition as:

\begin{corollary}\label{openpos3}   
Let  $A$ be an $n\times m$ matrix of rank $s$ and $Z$ an $m\times n$  influence specification. 
Then  $A$   is  injective over  $\mmK^g_{m,n}(Z)$  if and only if  the following two statements hold:
\begin{enumerate}[(i)]
\item  $A$ is injective over $\mmK^g_{m,n}[V]$ for some kinetic order $V$ with $Z(V)=Z$.
\item For all sets $I,J\subseteq \{1,\dots,n\}$ of cardinality $s$, if $A_{I,J}$ is non-singular then  $\det(V_{J,I})$  has the same fixed sign for all $V$ in $\Sigma(Z)$.
\end{enumerate}
\end{corollary}
Item $(ii)$ can be replaced by: $(ii)'$ $Z$ has a signed $A$-determinant.
In particular the proposition is true by choosing the kinetic order $V=V(Z)$.
The first condition  guarantees that all  coefficients in the polynomial expression of the determinant $\det(J_c(\widetilde{f}_{\kappa,V}))$  have the same sign or are zero, and that at least one coefficient is non-zero. The second condition ensures that this property is preserved for all kinetic orders $V$ in $\Sigma(Z)$. 

\begin{remark}\label{subset}
Using  Theorem~\ref{injecclose-pl} we have the following. If  $\mmK_0=\cup_{V\in {\bf V}}K^{g}_{m,n}[V]$ for some set ${\bf V}$ and the determinant $\det(J_{c}(\widetilde{f}_{\kappa,V}))$ does not vanish for all $(\kappa,V)\in\mmK_0$, then $A$ is injective over  $\mmK_0$, irrespectively whether $Z$ has a signed $A$-determinant or not. For example, if $\det(J_{c}(\widetilde{f}_{\kappa,V}))$ is not zero provided that $v_{j,1}>v_{j,2}$ for some $j$, then $A$ is injective over $\mmK_0=\{K^{g}_{m,n}[V]|\ v_{j,1}>v_{j,2},\ Z(V)=Z\}$. 
\end{remark}

The  next proposition  provides a characterization of injectivity of a matrix $A$ over  $\bigcup_{Z \vert Z_1\preceq Z\preceq Z_2}\mmK^{g}_{m,n}(Z)$.  

\medskip
\begin{proposition}\label{openpos}
Let  $A$ be an $n\times m$ matrix of rank $s$ and $Z_1\preceq Z_2$ two   $m\times n$  influence specifications. 
The following statements are equivalent:
\begin{enumerate}[(i)]
\item  $A$   is  injective over  $\bigcup_{Z \vert Z_1\preceq Z\preceq Z_2}\mmK^{g}_{m,n}(Z)$.
\item $Z_2$ has a signed A-determinant and $A$ is injective over $\mmK^{g}_{m,n}[V_1]$ for some $V_1$ in $\Sigma(Z_1)$.
\item  $A$   is  injective over  $\mmK^{g}_{m,n}(Z_1)$ and   $\mmK^{g}_{m,n}(Z_2)$.
 \end{enumerate}
\end{proposition}

 In particular the proposition is true for the kinetic order $V_1=V(Z_1)$.
The first part of (ii) guarantees that all  coefficients in the polynomial expression of the determinant $\det(J_c(\widetilde{f}_{\kappa,V_2}))$
have the same sign or are zero for all $V_2$ such that $Z(V_2)=Z_2$. Hence, by Lemma~\ref{support}(ii), this property is preserved for all kinetic orders $W$ with $Z(W)\preceq Z_2$. The second part of (ii)  ensures that at least one term is non-zero for all $W$  with $Z_1\preceq Z(W)$.

A natural choice for the smaller influence specification $Z_1$ is in many contexts the complex dependent influence specification $Z_{\mmC}$ given by the kinetic order $Y$.  Corollary~\ref{openpos3} and the discussion above  also imply the following corollary. 

\begin{corollary}\label{non-injective}
Let $A$ be the  $n\times m$ stoichiometric matrix. Assume that $A$ is not  injective over the set of mass-action kinetics $\mmK^{g}_{m,n}(Y)$ and that the determinant  $\det(J_c(\widetilde{f}_{\kappa,Y}))$ is not  identically zero. Then $A$ is not injective over $\mmK^{g}_{m,n}[V]$ for any kinetic order $V$ such that $Z_{\mmC}\preceq Z(V)$.
\end{corollary}

\begin{example}\label{ex:support}
Consider the network in Example \ref{futile} with the kinetic order $V$ introduced in Example \ref{ex:det2}. 
The  kinetic orders in $\Sigma(Z(V))$ are positively parameterized by the matrix
{\small $$Z_X=\left(  \begin{array}{cccc}  x_1 & x_2 & 0 & 0 \\ 0 & 0 & 0 & x_3 \\ 0 & 0 & 0 & x_4 \\ \sign(v)x_7 & 0 & x_5 & 0  \end{array} \right) $$}
The kinetic orders $W$ with $Z_{\mmC}\preceq Z(W)\preceq Z(V)$  also include the possibility $x_7=0$.   
The polynomial $p_{Z}$ corresponding to the determinant of $\widetilde{AZ_X}$  is
$$(x_1+x_2+x_3+x_4)x_5+x_2x_4+\sign(v)x_7 x_2.$$
If $v\geq 0$   then all  coefficients are positive.  In that case it follows from Proposition~\ref{openpos2} and Lemma~\ref{support}(i) that  the network is injective over  $\mmK^g_{m,n}(Z(V))$  and in particular over   $\mmK^g_{m,n}(Z_{\mmC})$ (corresponding to $v=0$).  If $v<0$ then the term  $\sign(v) x_7 x_2$ is negative while the rest are positive and the network is not injective over $\mmK^g_{m,n}(Z(V))$.  
If $v<0$ and the kinetic order $W$ fulfills $x_4>x_7$, then all terms are positive. It follows that the network is injective over $\mmK^{g}_{m,n}[W]$ (Proposition~\ref{det-criterion}, see also Remark~\ref{subset}) even though it is not injective over $\mmK^g_{m,n}(Z(V))$.
\end{example}

\begin{example}\label{futile:reactiondep}
Consider the network in Example \ref{futile} with the reaction dependent influence specification $Z_{\mmR}$. The  kinetic orders $V$ with $Z(V)=Z_{\mmR}$ are positively parameterized by
{\small $$Z_X = \left( \begin{array}{cccc} 
x_1 &  x_2 & 0& -x_3 
\\   -x_4  &  -x_5 & 0 & x_{6} 
 \\  -x_7 & 0  & -x_8 &  x_9 
 \\ 0 & -x_{10} & x_{11} & 0 
\end{array} \right)$$}
and $x_i>0$.
The monomials in the determinant of $\widetilde{AZ_X}$  have  positive coefficients. It follows that Proposition~\ref{openpos2} holds and that the network is injective over  $\mmK^{g}_{4,4}(Z_{\mmR})$. Example\,\ref{ex:support} showed that the network is injective over $\mmK^{g}_{4,4}(Z_{\mmC})$. Hence, it follows from Proposition~\ref{openpos} that the network is injective over $
\bigcup_{Z \vert Z_{\mmC}\preceq Z\preceq Z_{\mmR}}\mmK^{g}_{4,4}(Z)$.

The Mathematica code implementing the algorithm to decide whether the Example \ref{futile} is injective over $\mmK^{g}_{4,4}(Z_{\mmR})$ is shown in Figure~\ref{code}.
\end{example}

\begin{figure}[!t]
\centering
{\small \begin{minipage}[h]{0.9\textwidth}
\textcolor{gray}{Define the  stoichiometric matrix $A$: }
\begin{verbatim}
A = Transpose[{{-1,-1,0,1},{1,1,0,-1},{1,0,1,-1},{0,1,-1,0}}];
\end{verbatim}

\textcolor{gray}{Define the   matrix $Z_X$: }
\begin{verbatim}
ZX = {{x[1],x[2],0,-x[3]},{-x[4],-x[5],0,x[6]},{-x[7],0,-x[8],x[9]},
       {0,-x[10],x[11],0}};
{s,n,lengthx,Mtilde} = {MatrixRank[R],Length[R[[1]]],11,A.ZX}; \end{verbatim}

\textcolor{gray}{Find a reduced basis of $\im(A)^{\perp}$:}
\begin{verbatim}
If[s<n,
  F = RowReduce[NullSpace[A]]; 
  For[i=1,i<=Length[F],i++,
    Mtilde[[Flatten[Position[F[[i]],x_/;x!=0]][[1]]]]=F[[i]];   ];
];  \end{verbatim}
 
\textcolor{gray}{Compute the determinant of $\widetilde{AZ_X}$:}
\begin{verbatim} 
det = Expand[Determinant[Mtilde]];
\end{verbatim}
 
\textcolor{gray}{Check the signs of the coefficients:}
\begin{verbatim} 
Rules = {}; 
monomials = Flatten[MonomialList[det]];
For[i=1,i<=lengthx,i++,AppendTo[Rules,x[i]->1];];
sign = DeleteCases[DeleteDuplicates[Sign[monomials/.Rules]],0];
If[Length[sign] == 1, Print["The network IS injective"], 
          Print["The network is NOT injective"]; ];\end{verbatim}
\end{minipage}}
\caption{The algorithm, implemented in Mathematica, to decide whether Example \ref{futile} is injective over $\mmK^{g}_{4,4}(Z_{\mmR})$. } \label{code}
\end{figure}

\begin{remark} 
In general a network  will not be injective over the set of all power-law kinetics. 
In the case of  Example~\ref{ex:det}, the term $t=vk_1k_4c_1^vc_3$ is the only term depending on $v$, that is, the kinetic order. It changes sign with $v$, whereas none of the other non-zero terms do. Hence, the network cannot be injective over all power-law kinetics. 
\end{remark}

\begin{remark}
The decomposition \eqref{cauchy2} of $\det(\widetilde{AZ_X})$ is precisely the \emph{core determinant} of $A$, as defined in \cite[Lemma 3.7]{helton:determinant} for the complex dependent influence specification.
\end{remark}

\subsection{Injectivity of systems defined by submatrices of $A$}
The conditions  presented in the propositions in the previous sections relate to submatrices of  $A$.  
For  $J\subseteq\{1,\ldots,n\}$, let $A_{*,J}$ ($V_{J,*}$) be the restriction of $A$ ($V$) to the columns (rows) with indices in $J$.

\begin{theorem}\label{subnet}
Let $A$ be an $n\times m$ matrix of rank $s$, $V$ an $m\times n$ kinetic order and $J$ an index set of size $s$, such that $A_{*,J}$ has rank $s$. Assume that $A$ is injective over $\mmK^{g}_{m,n}[V]$.
Then either 
\begin{enumerate}[(i)]
\item $A_{*,J}$ is injective over $\mmK^{g}_{s,n}[V_{J,*}]$ and all steady states are non-degenerate, or 
\item $A_{*,J}$ has only degenerate steady states.
\end{enumerate}
\end{theorem}

The theorem relates to  \cite[Cor.~8.1, Cor.~8.2]{Feliu-inj}, where injectivity of  a  network is studied relatively to  injectivity of the network augmented with the ``missing" outflow reactions.  Also, the theorem relates to Joshi and Shiu  \cite{joshi-shiu-II}. They consider a network obtained by restricting a larger network such that the stoichiometric dimension is maintained. If the smaller network has multiple steady states then so does the larger.
 Theorem~\ref{subnet} cannot be used to draw the same conclusion since non-injectivity does not imply that there are multiple steady states.

\begin{example}\label{ex:support2}  
 Consider Example~\ref{pop-growth}. The system has stoichiometric matrix $A=(1,-1)$, rank one and kinetic order $V=(1,2)$. There are two possible choices of index set, $J=\{1\}$ and $J=\{2\}$. The matrices $A_{*,i}$ are  injective over  $\mmK^g_{1,1}[V_{i,*}]$, $i=1,2$, but $A$ is not injective over $\mmK^g_{1,1}[V]$. 
\end{example}

\section{Injectivity for strictly monotonic kinetics}
\label{general-kinetics}

In this section we extend the results on injectivity for power-law kinetics in Sections \ref{sec:powerlawinj} and \ref{sec:int} to cover $Z$-injectivity of a matrix $A$ over the set of strictly monotonic kinetics  $\mmK_{m,n}(Z)$. The following theorem is the main theorem of this article.

\begin{theorem}\label{equiv}
Let  $A$ be an $n\times m$ matrix  and $Z$ an $m\times n$  influence specification.  The following three statements are equivalent:
\begin{enumerate}[(i)]
\item $A$ is $Z$-injective  over  $\mmK_{m,n}(Z)$.
\item $A$ is $Z$-injective  over $\mmK^{d}_{m,n}(Z)$.
\item $A$ is injective over  $\mmK^{g}_{m,n}(Z)$.
\item $Z$ is A-SNS.
\end{enumerate}
\end{theorem}

 The theorem implies that  for a matrix to be $Z$-injective over $\mmK_{m,n}(Z)$  it is sufficient to be injective over  $\mmK_0$ with $\mmK^{g}_{m,n}(Z)\subseteq \mmK_0\subseteq \mmK_{m,n}(Z)$. In  \cite[Prop.~5.2]{Feliu-inj} it is shown that injectivity and $Z$-injectivity are equivalent notions for mass-action kinetics. Theorem~\ref{equiv} also implies that  \cite[Prop.~5.2]{Feliu-inj} holds  generally, namely that $A$ is $Z$-injective  over  $\mmK_{m,n}(Z)$ if and only if $A$ is injective  over  $\mmK_{m,n}(Z)$.

\begin{theorem}\label{injecclose} Let  $A$ be an $n\times m$ matrix. Then the following are equivalent:
\begin{enumerate}[(i)]
\item  $\ker(J_{c}(f_{K}))\cap \im(A)=  \{0\}$ for all $c\in \R^n_+$ and $K\in\mmK^{d}_{m,n}(Z)$.
\item  $\ker(J_{c}(f_{\kappa,V}))\cap \im(A)=  \{0\}$ for all $c\in \R^n_+$ and $(\kappa,V)\in\mmK^{g}_{m,n}(Z)$.
\end{enumerate}
If either of the two statements holds then $A$ is $Z$-injective over $\mmK_{m,n}(Z)$.
\end{theorem}

Recall from \eqref{jaccomp} that  $\ker(J_{c}(f_{K}))\cap \im(A)=\{0\}$  if and only if $\det(J_{c}(\widetilde{f}_{K}))\neq 0$.

\begin{corollary}
If $A$ is injective over $\mmK^{d}_{m,n}(Z)$ then $A$ cannot have positive  degenerate steady states.
 \end{corollary}

The corollary follows immediately from \eqref{jaccomp}. The same result holds for weakly monotonic kinetics (\cite{shinar-conc}, see Definition~\ref{def-mono}) and for  mass-action kinetics \cite{Feliu-inj}.

\begin{example} 
According to Example~\ref{futile:reactiondep}, the network in Example~\ref{futile} is injective over $\bigcup_{Z \vert Z_{\mmC}\preceq Z\preceq Z_{\mmR}}\mmK^{g}_{m,n}(Z)$. It follows from Theorem~\ref{equiv} 
 that  the network is $Z$-injective over $\mmK_{m,n}(Z)$ for any $Z$ such that $Z_{\mmC}\preceq Z\preceq Z_{\mmR}$.
\end{example}

 \begin{remark}
The criterion that $Z$ is $A$-SNS is computationally efficient and can easily be implemented using symbolic software  (see also Figure~\ref{code}). It requires calculation of the matrix $\widetilde{AZ_X}$,  its determinant and the expansion of the determinant. The complexity of the latter depends on the number of species influencing a reaction as well as the size of the matrix, while the former  depends on the size of the matrix only. The criterion provided in Corollary~\ref{polyterms} is computationally more demanding as it requires investigating all minors of a certain size. 
\end{remark}

\section{Graphical representation of the criterion}\label{graphical_sec} 
We have shown that a matrix $A$ is injective over the set of kinetics strictly monotonic with respect to an influence specification $Z$, if and only if  $Z$ is $A$-SNS. The property of being $A$-SNS relies on  computing the symbolic determinant of $\widetilde{AZ_X}$.
Since visual inspection is often more appealing than computation, injectivity-based criteria to  preclude multistationarity have been interpreted in graph-theoretical terms \cite{craciun-feinbergII,banaji-craciun1,banaji-craciun2,Shinar-conc2}. Generally, the outcome does not provide a full characterization of injectivity but only a sufficient graphical condition that guarantees  injectivity of $A$.
In this section we show  that the decomposition \eqref{cauchy2} can be  interpreted directly in terms of circuits in the DSR-graph.

The procedure basically relies on a variant of the DSR-graph (directed-species-reaction-graph) introduced in \cite{banaji-craciun2}, using  the matrices $A$ and $Z_X$. To fix the notation, consider two sets  $\mmS=\{S_1,\dots,S_n\}$ (``species'') and  $\mmR=\{r_1,\dots,r_m\}$ (``reactions''). 
The DSR-graph, $G_{A,Z}$, associated to $(A,Z_X)$  is defined in the following way. The set of nodes of the graph is $\mmS\cup \mmR$ and hence there are $n+m$ nodes. 
There is a directed edge from a species $S_i$ to a reaction $r_j$ if and only if $z_{j,i}\neq 0$.
 This edge is assigned a symbolic label $e_{j,i}:=z_{j,i}\cdot x_{j,i}$.
There is a directed edge from a reaction $r_j$ to a species $S_i$ with label $a_{i,j}$ if and only if $a_{i,j}  \neq 0$.

 A circuit in a graph $G$ is a sequence of distinct nodes $i_1,\ldots,i_l$ such that there is a directed edge from $i_k$ to $i_{k+1}$ for all $k\leq l-1$ and one from  $i_l$ to $i_1$.  A circuit must involve at least one edge.   In this specific setting,  any circuit involves an even number of edges. The label of a circuit is the product of the labels of the edges in the circuit. Two circuits are disjoint if they do not have any common nodes.
 A circuit has sign $(-1)$ if the number of species nodes in the circuit is even.

A $k$-nucleus  is a collection of disjoint circuits which involves $k$  nodes. 
 The label $l(D)$ of a  $k$-nucleus $D$ is the product of the labels of the edges in the nucleus.
The sign of a $k$-nucleus is $(-1)^q$, where $q$ is the number of circuits with even number of species nodes.
That is, if  $D=C_1\cup \dots \cup C_a$ is  a disjoint union of circuits then
$$\sign(D)l(D) = \prod_{i=1}^a \sign(C_i) l(C_i).  $$

\begin{proposition}\label{graphical}
Fix two sets $I,J\subseteq \{1,\dots,n\}$ of cardinality $s$ and consider the submatrices $A_{I,J}$ and $(Z_X)_{J,I}$. 
Let $D_s(I,J)$ be the set of $2s$-nuclei of $G_{A,Z}$ with nodes $S_i$ for $i\in I$ and $r_j$ for $j\in J$. 
Then
$$\det(A_{I,J})\det( (Z_X)_{J,I}) = \sum_{D \in D_s(I,J)} \sign(D)l(D).  $$

\end{proposition}

It follows from the proposition that injectivity of  $A$ with respect to $Z$ can be decided from the DSR-graph as follows:
\begin{enumerate}[(1)]
\item Classify all  circuits $C$ of the graph $G_{A,Z}$ according to the number of species nodes that are involved (up to $s$). Assign  $\sign(C)l(C)$ to each of them.
\item Consider all products of circuit labels with sign  for which the number of species nodes adds up to $s$.
\item Each label has the form $\lambda x_{j_1,i_1}\cdot \dots \cdot x_{j_s,i_s}$ for some scalar $\lambda$. Keep only the monomials for which there is no repetition among $i_1,\dots,i_s$ and among $j_1,\dots,j_s$.
\item Add the terms with the same variables.
\end{enumerate}

According to Proposition~\ref{graphical}, the terms obtained after step (4) are exactly the terms in the decomposition of the determinant \eqref{cauchy2}.
Using this representation non-necessary conditions for injectivity might be developed by relaxing the information initially encoded in the DSR-graph, in the spirit of \cite{banaji-craciun1,banaji-craciun2,craciun-feinbergII,Shinar-conc2}.

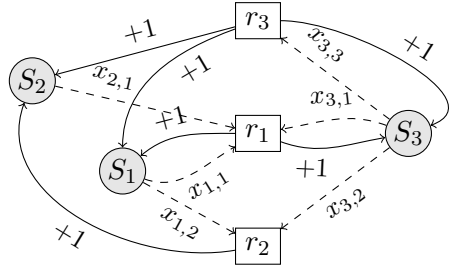
\begin{figure}[t]
\centering
\begin{tikzpicture}
\node[draw,circle,fill=gray!20!white,inner sep=1pt] (A) at (-1.8,-0.5) {$S_1$};
\node[draw,circle,fill=gray!20!white,inner sep=1pt] (B) at (-3,0.7) {$S_2$};
\node[draw,circle,fill=gray!20!white,inner sep=1pt] (C) at (2,0) {$S_3$};
\node[draw] (r1) at (0,0) {$r_1$};
\node[draw] (r2) at (0,-1.5) {$r_2$};
\node[draw] (r3) at (0,1.5) {$r_3$};
\draw[->,dashed] (A)  to node[below,sloped] {\small $x_{1,2}$} (r2);
\draw[->,dashed] (A)  to[out=340,in=210] node[below,pos=0.6,sloped] {\small $x_{1,1}$} (r1);
\draw[->,dashed] (B)  to node[above,pos=0.3,sloped] {\small $x_{2,1}$} (r1);
\draw[->,dashed] (C)  to node[below,sloped] {\small $x_{3,2}$} (r2);
\draw[->,dashed] (C)  to[out=160,in=10] node[above,sloped] {\small $x_{3,1}$} (r1);
\draw[->,dashed] (C)  to node[above,pos=0.7,sloped] {\small $x_{3,3}$} (r3);
\draw[->] (r2)  to[out=190,in=250] node[below,sloped] {\small $+1$} (B);
\draw[->] (r3)  to  node[above,sloped] {\small $+1$} (B);
\draw[->] (r3)  to[out=200,in=90]  node[below,pos=0.3,sloped] {\small $+1$} (A);
\draw[->] (r1)  to[out=180,in=40]  node[above,pos=0.6,sloped] {\small $+1$} (A);
\draw[->] (r3)  to[out=0,in=20]  node[above,sloped] {\small $+1$} (C);
\draw[->] (r1)  to[out=340,in=190]  node[below,pos=0.3,sloped] {\small $+1$} (C);
\end{tikzpicture}
\caption{DSR-graph of Example~\ref{ex:dsr}.}\label{fig2}
\end{figure}

\begin{example}\label{ex:dsr}
Consider the network   defined by the reactions $r_1\colon S_1+S_2+S_3\rightarrow 2S_1+S_2+2S_3$, $r_2\colon S_1+S_3\rightarrow S_1+S_2+S_3$, $r_3\colon S_3\rightarrow S_1+S_2+2S_3$ and  the complex dependent influence specification, $Z=Z_{\mmC}$. In this case $s=2$. The DSR-graph of this network is shown in Fig. \ref{fig2}.
We enumerate all circuits of the graph and classify them according to the number of species nodes:  

\medskip
\begin{center}
\begin{tabular}{c|c|c}
$\#$ Species nodes & 1 & 2  \\ \hline
Label & $x_{1,1}$, $x_{3,1}$, $x_{3,3}$ & $-x_{2,1}x_{3,2}$, $-x_{2,1}x_{3,3}$, $-x_{1,1}x_{3,3}$, $-x_{1,2}x_{2,1}$ 
\end{tabular}
\end{center}

\medskip
We take the products of circuit labels for which the number of species nodes adds to $s=2$ and avoid index repetition:
$$-x_{2,1}x_{3,2},\quad -x_{2,1}x_{3,3}, \quad -x_{1,1}x_{3,3}, \quad -x_{1,2}x_{2,1},\quad \text{ and } \quad x_{1,1}x_{3,3}$$
 (the latter is obtained by combining two circuits with one species node). Note that the circuit with label $x_{3,1}$ cannot be combined with any other circuit because there would be index repetitions.
We group the monomials together whereby the coefficient of $x_{1,1}x_{3,3}$ becomes zero. The remaining three monomials are the terms of the polynomial $p_Z(X)$.
Since the terms have  the same sign, the network is injective.
\end{example}

\section{Extensions to other types of  influence specifications}\label{sec:shinar}

Shinar and Feinberg \cite{shinar-conc}  introduce the term ``weakly monotonic kinetics", which in some sense imposes a weaker requirement on the kinetics  than the term ``strictly monotonic kinetics" introduced here (Definition~\ref{def-mono0} and Lemma~\ref{lemma-mono0}). In this section we assume  that 
we have a network $\mmN=(\mmS,\mmC,\mmR)$ and an influence specification $Z$ such that $\supp(y_j)\subseteq z^+_{j}\cup z^-_{j}$. That is, the enhancers and inhibitors of a reaction include all species involved in the reactant  complex $y_j$. In chemical reaction theory it is typically required that the species in the reactant complex have positive influence on the reaction \cite{gunawardena-notes,shinar-conc} and not  negative or neutral, as in Definition~\ref{infl-spec}. However, relaxation of this assumption is found, for example in S-systems theory. 

\begin{definition}[\cite{shinar-conc}]
\label{def-mono}
A  kinetics  $K$ for a network $\mmN=(\mmS,\mmC,\mmR)$  is \emph{weakly monotonic with respect to  an influence specification} $Z$ if, for each pair of non-overlapping concentration vectors $a,b\in\Omega_K$, the following implications hold for all   $j$: 
\begin{enumerate}[(i)]
\item if $K_{j}(a)>K_{j}(b)$ then  $\sign(a_i-b_i)=z_{j,i}\not=0$ for some $i$.
\item if $K_{j}(a)=K_{j}(b)$ then $a_i=b_i$ for all  $i\in\supp(y_j)$, or  $\sign(a_i-b_i)=z_{j,i}\neq 0$ and $\sign(a_{\ell}-b_{\ell})=-z_{j,\ell}\not=0$ for some distinct $i,\ell$.
\end{enumerate}
Let $\mmK^w_{m,n}(Z)$ denote the set of   kinetics that are weakly monotonic with respect to $Z$. 
\end{definition}

Using the characterization of strictly monotonic kinetics provided in Lemma~\ref{lemma-mono0}, we find that the two definitions differ in  (ii), where it is required that $a_i=b_i$ for all $i\in z^+_{j}\cup z_{j}^-$ and not just for $i\in\supp(y_j)$. In this sense, Definition~\ref{def-mono} imposes a weaker requirement on the kinetics than Definition~\ref{def-mono0} and 
$$\mmK_{m,n}(Z) \subset \mmK^w_{m,n}(Z). $$ 
 In \cite{shinar-conc}, an  influence specification $Z$ must fulfill  $z_{j,i}=+$ for $i\in\supp(y_j)$. 
 Definition \ref{def-mono0} stipulates that all species   play an equal role in the definition, whereas Definition  \ref{def-mono} singles out the species in the reactant complex to have special importance. Below we show that our determinant criterion also  applies to the broader definition of influence specification.

The determinant criterion  in Theorem~\ref{equiv}(iv)  can be adapted to derive a determinant criterion for a network $\mmN$  to be $Z$-injective over $\mmK^w_{m,n}(Z)$. We first note that the influence specification $\widetilde{Z}$ given by $\widetilde{z}_{j,i}=z_{j,i}$ for $i\in\supp(y_j)$ and zero otherwise is a minimal element among all influence specifications $Z'$ for $\mmN$  fulfilling $\supp(y_j)\subseteq z^+_{j}\cup z^-_{j}$  and $Z'\preceq Z$.
Also note that if $a,b$ are $Z$-non-overlapping then they are $Z'$-non-overlapping for all influence specifications $Z'$ such that $Z'\preceq Z$. Hence, we have that
\begin{equation}\label{w-implies}
Z'\preceq Z \quad \text{ implies } \quad \mmK^w_{m,n}(Z')\subseteq \mmK^w_{m,n}(Z).
\end{equation}
 
\begin{lemma}\label{gma-influence} 
If $K\in\mmK^w_{m,n}(Z)$ is a power-law kinetics then  there is $\widetilde{Z}\preceq Z'\preceq Z$ such that $K\in\mmK^g_{m,n}(Z')$. That is, 
$$\mmK^w_{m,n}(Z) \cap \mmK^g_{m,n}=\bigcup_{Z' \vert\widetilde{Z}\preceq Z'\preceq Z}\mmK^g_{m,n}(Z').$$
\end{lemma}

\begin{theorem}\label{equivw} 
 Let $\mmN=(\mmS,\mmC,\mmR)$ be a network with influence specification $Z$.  The following statements are equivalent:
\begin{enumerate}[(i)]
\item $\mmN$ is $Z$-injective  over  $\mmK^w_{m,n}(Z)$.
\item $\mmN$ is $Z$-injective  over  $\bigcup_{Z' \vert\widetilde{Z}\preceq Z'\preceq Z}\mmK_{m,n}(Z')$.
\item $\mmN$ is injective over  $\bigcup_{Z' \vert\widetilde{Z}\preceq Z'\preceq Z}\mmK^g_{m,n}(Z')$.
\end{enumerate}
\end{theorem}

\medskip
Together with Proposition~\ref{openpos} we can derive a determinant criterion for a network $\mmN$ to be $Z$-injective  over  $\mmK^w_{m,n}(Z)$. Further, it is straightforward to derive statements similar to those in Section \ref{general-kinetics} for $\mmK^w_{m,n}(Z)$.

\begin{remark}
Shinar and Feinberg \cite{shinar-conc}  introduce the concepts of a \emph{concordant network} and of a network being \emph{concordant with respect to an influence specification $Z$}. To be concordant depends on the kinetics associated with the network only through the influence specification $Z$.   
They show that to be concordant with respect to $Z$ is equivalent to be injective over $\mmK^w_{m,n}(Z)$. Theorem~\ref{equivw} and Theorem~\ref{equiv}  provide an equivalent characterization in terms of  the influence specification through properties of the matrices $Z$ and $A$, and bring out an explicit relationship to the set of power-law kinetics. 
\end{remark}

 \begin{remark}
For Windows-based platforms, the CRN Toolbox \cite{crnttoolbox} provides  a test for injectivity over $\mmK^w_{m,n}(Z)$ for the influences introduced in \cite{shinar-conc}. 
\end{remark}

\section{The $P$-matrix property}\label{sec:banaji}
\label{p-matrix}
In \cite{banaji-donnell,banaji-craciun2}, an injectivity-related criterion is given to preclude the existence of multiple steady states in  $\mmC^1$  dynamical systems admitting a decomposition of the form $f(c)=AK(c)$.
The kinetics $K$ is required to be \emph{non-autocatalytic} (NAC), which is a condition that also involves the form of $A$. We introduce it in terms of a corresponding network $\mmN=(\mmS,\mmC,\mmR)$.
 In our terminology, NAC implies that (a) no species appear both in the reactant and the product complex of a reaction and (b) the influence specification fulfills $Z\preceq Z_{\mmR}$.  
   This class of dynamical systems includes the  differentiable kinetics $\mmK^d_{m,n}(Z)$ (assuming further that the rate functions are $\mmC^1$ in the domain of differentiability). In particular,  multistationarity in networks with power-law kinetics can be precluded using the criterion. 
 
The focus is on conditions for injectivity of the ``open network'' which is the  network obtained by adding the outflow reactions $S\rightarrow 0$, $S\in\mmS$, to the network, unless they are already present. The species $S$ is required to be the only species with non-zero (positive) influence on the reaction. If the reaction is already in the network, it must fulfill this requirement too.

If the open network can be constructed and is injective, then the initial network cannot have multiple non-degenerate steady states  \cite{banaji-craciun2}. Therefore, multistationarity can be precluded in an arbitrary  network provided it can be precluded in the corresponding open network. 
However  networks exist that are injective but for which the corresponding open network  is not injective (see Example~\ref{counterexample}).

Injectivity of the open network  follows from the results of Gale and Nikaid{\^o} \cite{Gale:1965p474} after determining that the Jacobian of the system associated to the open network is a $P$-matrix, or, equivalently, that the Jacobian of the system associated to the initial network is a $P_0$-matrix  \cite{banaji-craciun2}. A square matrix is said to be a $P$-matrix if all principal minors of the matrix are positive. If the principal minors are non-negative then the matrix is said to be a $P_0$-matrix. 
We proceed to discuss the relationship between  the $P_0$-matrix property and our criteria.

The next proposition is established in \cite{banaji-craciun2}  (stated using our terminology).
\begin{proposition}[\cite{banaji-craciun2}, Lemma 3.5]\label{banaji}
Let $A$ be a stoichiometric matrix, $Z$ an influence specification and  $V$ a kinetic order such that $Z(V)\preceq Z$. If 
\begin{itemize}
\item[$(*)$]$(-1)^{\ell}\det(A_{I,J})\det(V_{J,I})\geq 0$ for all sets $I,J\subseteq\{1,\ldots,n\}$ of cardinality $n-\ell$, for all $\ell$, 
\end{itemize}
then $-J_c(f_{\kappa,V})$  (minus the Jacobian) is a $P_0$-matrix.
\end{proposition}

If we require the kinetics  to be NAC   then $Z(V)\preceq Z_{\mmR}$ and the corresponding entries in  $-V$ and $A$ have the same sign. 
The determinant criterion in Proposition~\ref{det-criterion}(ii) (without the requirement that one product is non-zero) is implied by criterion $(*)$. Also if $(*)$ holds for all kinetic orders $V$ such that $Z(V)=Z$ then  the criterion in Corollary~\ref{openpos3}(ii)
is implied by criterion $(*)$. However,  as Example~\ref{counterexample} below shows, the opposite is not true. Hence, our criterion is  weaker than  criterion $(*)$.
The additional requirement in Proposition~\ref{det-criterion} and Corollary~\ref{openpos3}(ii)
 that one term is non-zero is necessary (and sufficient) to guarantee that $A$ is injective. 
In fact, if all terms are zero then all steady states are degenerate (see \cite{Feliu-inj} for a discussion that relates injectivity of the open network to that of the initial network for mass-action kinetics).

In \cite{banaji-craciun2}, the authors further provide a graphical condition  on the DSR-graph that implies that the Jacobian is a $P$-matrix and consequently that the network is injective.

\begin{example}\label{counterexample}
Consider the  network $\mmN$ defined by the set of reactions $\mmR=\{S_1+S_2+S_3\rightarrow 2S_1+S_2+2S_3, S_1+S_3\rightarrow S_1+S_2+S_3, C\rightarrow S_1+S_2+3S_3\}$ and  complex dependent influence specification, $Z=Z_{\mmC}$. For a kinetic order $V$ with $Z(V)=Z$, let
$$
V = \left( \begin{array}{ccc} 
v_{1,1} & v_{2,1} & v_{3,1} \\  v_{1,2}  & 0 & v_{3,2} \\ 0 & 0 & v_{3,3}   
\end{array} \right),\quad \quad
A = \left( \begin{array}{ccc} 
1 & 0 & 1\\  0  & 1 & 1 \\ 1 & 0 & 2
\end{array} \right).$$
The stoichiometric space has maximal dimension and
 $$\det(AV)=-v_{1,2}v_{2,1}v_{3,3}.$$ 
 Consequently, $\mmN$ (or $A$) is injective over  $\mmK^g_{m,n}(Z)$ (Proposition~\ref{openpos2}), hence also over $\mmK_{m,n}(Z)$ (Theorem~\ref{equiv}). However, the product of the minors obtained by removing the second row and column is positive, while that obtained by removing the third row and column is negative. 
Hence, $-J_c(f_{\kappa,V})$  is not  a $P_0$-matrix. In fact, the open network associated to $\mmN$ with influence specification $Z$ is not injective.
\end{example}

\begin{example}\label{counterexample2}
Consider Example~\ref{ex:dsr}, which is similar Example~\ref{counterexample}, but with  the coefficient of $S_3$ changed in the last reaction. The stoichiometric space has dimension two. All non-zero  products involving $2\times 2$ matrices in criterion $(*)$  have negative sign and, hence, $(*)$ is not fulfilled.   However, since all  products have the same sign (as we saw in Example~\ref{ex:dsr}),  then the network is injective over $\mmK_{m,n}(Z)$.
\end{example}

\section{The interaction graph} \label{sec:soule} 

Conditions for the preclusion of multistationarity have also been given for generic  dynamical systems described by ordinary differential equations and we will here review one condition due to  Kaufman, Soul\'{e} and Thomas \cite{kaufman,soule} that closely relates to our work.  This condition is also based on the Jacobian of the system and takes the form of a graphical condition.  Specifically, we will interpret a result of \cite{kaufman} in terms of our framework and show, by example, that our criterion might decide on injectivity when  the criterion in \cite{kaufman} fails. For this we need some preliminaries.
 
 Let a dynamical system $\dot{c}=F(c)=(F_1(c),\ldots,F_n(c))$ be given such that $c=(c_1,\dots,c_n)\in \Omega\subseteq \R^n$, where $\Omega$ is a product of open intervals of $\R$ and $F_j$ is differentiable in the interior of $\Omega$. The interaction graph $G(c)$ at $c$ is the labeled directed graph with node set $\{1,\ldots,n\}$ and labels in  the set $\{-,+\}$ such that there is  an edge from node $i$ to node $j$ if $\partial F_j(c)/\partial c_i\neq 0$. The edge has label given by $\sign(\partial F_j(c)/\partial c_i)$.
Denote by $\widehat{G}(c)$  the $n\times n$ sign matrix with $(j,i)$th entry  $\sign(\partial F_j(c)/\partial c_i)$. It encodes the same information as $\widehat{G}(c)$.  We use the definition of a circuit and a $k$-nucleus given in Section~\ref{graphical_sec}. However, in this section, 
the sign of a circuit is the product of the labels of the edges in the circuit. The sign of a   $k$-nucleus is $(-1)^{p+1}$ where $p$ is the number of  circuits in the $k$-nucleus with sign equal to $+$, i.e. positive \cite{soule}.  A $k$-nucleus is variable if one edge in one of the circuits does not have constant sign in $c$. 

In \cite[Th.~2]{kaufman} (see below), a mild regularity condition is imposed  on the function $F$. To keep the presentation clear, the reader is referred to the original paper for its description. We refer  to it  as condition (C). 

\begin{theorem}[\cite{kaufman}]\label{kaufman}
Assume that the system $\dot{c}=F(c)$ has two non-degenerate steady states and that $F$ fulfills condition (C). Then one of the following statements is true:
\begin{enumerate}[(i)]
\item There exists $c\in\Omega$ such that $G(c)$ has two $n$-nuclei of different sign.
\item There is a variable $n$-nucleus.
\end{enumerate}
\end{theorem}

The existence of non-degenerate steady states implies that the Jacobian of $F$ is non-singular for all $c$ and  that a $n$-nucleus exists for some $c$.
If neither (i) nor (ii) above are fulfilled then the system cannot have multiple non-degenerate steady states. 
There is no a priori restriction to pairs of \emph{positive} steady states (unless $\Omega\subseteq\R^n_+$). 
If (ii) is not fulfilled then preclusion of multiple steady states must follow from the failure of (i). Therefore,  in relating our work to \cite{kaufman},  we assume that $G(c)$ does not depend on $c$, that is $G(c)=G$ for all $c$ (and (ii) is not fulfilled by hypothesis).

Theorem~\ref{kaufman} can be rephrased in  our terminology as a  statement about  the preclusion of multiple non-degenerate steady states  in any dynamical system with interaction graph $G$.  In particular, we show that failure of condition (i) is equivalent to  non-injectivitity of  certain matrices $A$ over the class of differentiable kinetics. 

In what follows we assume that the first $s$ rows of $\widehat{G}$ are non-zero and that the last $n-s$ rows are identically zero. This can always be obtained by permuting the order of the variables $c_1,\dots,c_n$. If $G$ is the (constant) interaction graph associated to a  dynamical system as above, then a zero row of $\widehat{G}$ corresponds to a constant $F_j$. If a zero row exists, then the Jacobian of $F$ is necessarily singular and  Theorem~\ref{kaufman} cannot be applied to preclude multistationarity. 
In general, there might be many decompositions $A, Z$  of the system, such that the sign pattern of $AZ$ is $\widehat{G}$.  Theorem~\ref{kaufman} does not distinguish between these. In order to relate the theorem to our setting, we introduce a family of decompositions $A$ and $Z$. 

\begin{definition}\label{netgraph} 
 Let $\widehat{G}=(g_{i,j})_{i,j}$ be an $n\times n$ sign matrix with non-zero rows $1,\dots,s$. 
For each $j=1,\dots,s$, choose a set $H_j\subseteq \{1,\ldots,s\}$.  We associate with these sets an $n\times (2s)$ stoichiometric matrix $A=(a_{i,j})_{i,j}$ and an  $(2s)\times n$ influence specification $Z=(z_{j,i})_{j,i}$  by 
\begin{align*}
a_{j,j}  &= 1,   & z_{j,i} &=  g_{j,i},  && \hspace{-1cm} i\in H_j \\
a_{j,j+s} &=-1,  &  z_{j+s,i} &=  -g_{j,i}, &&  \hspace{-1cm} i\in \{1,\ldots,s\}\setminus H_j
\end{align*}
for $i=1,\ldots,n$, $j=1,\ldots,s$, and zero otherwise.  
\end{definition}

In other words, the  top $s\times (2s)$ block of  $A$ is composed of two diagonal matrices adjacent to each other and the bottom $n-s$ rows are zero. $Z$ consists of  the first $s$ rows of $\widehat{G}$ duplicated, but with a change of sign in some entries. 
The definition can be casted in terms of a reaction network with species $\{S_1,\ldots,S_n\}$ and set of reactions  $\mmR=\{ 0\rightarrow S_j, S_j\rightarrow 0 \vert \ j=1,\ldots,s\}$. This network has stoichiometric matrix $A$ with rank $s$.
Clearly, the sign pattern of $AZ$ is $\widehat{G}$ by construction.
Furthermore, for any $K\in \mmK^d_{2s,n}(Z)$, the system $f_K(c)=A K(c)$ has interaction graph $G$.

\begin{theorem}\label{soule}
Let $G$ be an interaction graph and define $A, Z$, and $s$ as in Definition~\ref{netgraph}. 
The following two statements are equivalent:
\begin{enumerate}[(i)]
\item $A$ is $Z$-injective over $\mmK^d_{2s,n}(Z)$.
\item $G$ has at least one $s$-nucleus and all $s$-nuclei of $G$ have the same sign.
\end{enumerate}
\end{theorem}

Note that (ii) is independent of the choice of $H_j$ in Definition~\ref{netgraph}, and hence  statement (i) is also independent of the choice.

Assume that $s=n$ and $A$ is $Z$-injective with $A, Z$ chosen as in Definition~\ref{netgraph}. Then any dynamical system $\dot{c}=F(c)$ with associated interaction graph $G$ cannot have multiple non-degenerate positive steady states. If $F$ fulfills condition (C), then by Theorem~\ref{soule} and Theorem~\ref{kaufman}(i), multiple non-degenerate steady states are precluded (not only positive steady states).
This is, in particular, true if we choose $H_j=\{1,\ldots,s\}$ for all $j$. In this case, we might choose   $A=I_n$ (the $n\times n$ identity matrix) and $Z=\widehat{G}$, as the bottom half of $Z$ in Definition~\ref{netgraph} is identically zero.
These results lead to the following corollary  (using Proposition~\ref{openpos2}), which was proven in \cite{gouze}, following a more direct route.

\begin{corollary}\label{graphsns}
Let $\widehat{G}$ be an $n\times n$ $I_n$-SNS matrix. Then, any dynamical system $\dot{c}=F(c)$ in $\R^n$  that fulfills condition (C) and  has constant sign matrix $\widehat{G}$  cannot have multiple non-degenerate steady states.
\end{corollary}

 We finish the section with two examples that   illustrate the relationship between the criteria.

\begin{example}
Consider Example~\ref{gene-regu}. Karlebach and Shamir \cite{karlebach} model the gene network depicted in Figure~\ref{Xgenes} as
$$\dot{c}_1=\frac{\alpha_1}{1+\beta_1 c_3}-\delta_1 c_1, \quad \dot{c}_2=\frac{\alpha_2 c_1}{1+\beta_2 c_1}-\delta_2 c_2, \quad
\dot{c}_3=\frac{\alpha_3 c_1c_2}{(1+\beta_3 c_1)(1+\beta_4 c_2)}-\delta_3 c_3,$$
for positive parameters, $\alpha_i,\beta_i,\delta_i$. Here degradation of each gene $G_1,G_2,G_3$ is incorporated.
The interaction graph $G$ associated with the system is constant and
$$\widehat{G} = \left(\begin{array}{rrr} - & + & + \\ 0 & - & + \\ - & 0 & -  \end{array}\right).$$
This matrix is $I_n$-SNS ($n=3$) and, therefore, the system cannot have multiple non-degenerate steady states.

Using the stoichiometric matrix $A=I_n$ effectively corresponds to analyzing injectivity of a network with only inflow reactions ($0\to G_i$) and influence specification given by $\widehat{G}$. Since, in our terminology, a kinetics must be positive, the entry $-$ in position $(1,1)$ of $\widehat{G}$ corresponds to a decreasing kinetics in $c_1$ for the inflow reaction $0\rightarrow G_1$. This is different from the system we started from, which had a negative summand ($-\delta_1c_1$). Alternatively, each $\dot{c}_j$ might be separated into two components, one representing the reaction $0\to G_i$, the other $G_i\to 0$.
\end{example}

As we can deduce from the results above, preclusion of multistationarity by the methods of \cite{kaufman}, is essentially preclusion of multistationarity in networks in which only inflow reactions are considered. Knowledge about the underlying network structure allows us to  preclude multistationarity for a bigger class of dynamical systems. We do not only ``see'' the signs of the entries of the Jacobian, but also the terms that contribute to the signs.  This is illustrated in the following simple example.

\begin{example}\label{soule-ex}
Consider the sign matrix $\widehat{G}$ with entries $g_{1,2}=g_{2,1}=+$ and $g_{1,1}=g_{2,2}=-$.
This matrix is not $I_n$-SNS and hence multistationarity cannot be excluded from Theorem~\ref{kaufman}.
Consider now a network with  reactions  $r_1\colon S_1\rightarrow S_2$ and $r_2\colon S_2\rightarrow 0$, and stoichiometric matrix $A$ given by the reactions $r_1,r_2$. The dimension of the stoichiometric subspace is $n=2$.
For any kinetics $K\in\mmK_{2,2}$, the ODE system associated with the network is of the from
$$\dot{c}_1 =- K_1(c) \qquad \dot{c}_2 =  K_1(c) - K_2(c). $$
Consider the influence specification $Z$ with non-zero terms:  $z_{1,1}=z_{2,2}=+$ and $z_{1,2}=-$. If $K\in \mmK^d_{2,2}(Z)$, then the  Jacobian of the species formation rate function has interaction graph $G$. Further, the only set of reactions of cardinality  $s=n=2$ is $\{r_1,r_2\}$. By calculation,   $Z$ is $A$-SNS and  it follows from Proposition~\ref{openpos2} and Theorem~\ref{equiv}  that the network is $Z$-injective over $K^d_{2.2}(Z)$ and multistationarity cannot occur.
\end{example}

\section{Hill-type kinetics and injectivity}\label{sec:hill}

Let $A$ be an $n\times m$ stoichiometric matrix and $Z$ an $m\times n$ influence specification. The key to the statements in Section~\ref{general-kinetics} is that whenever there are two non-overlapping stoichiometrically compatible concentration vectors $a,b$ and $f_K(a)=f_K(b)$ for some kinetics $K\in\mmK_{m,n}(Z)$, then we can find two positive stoichiometrically compatible concentration vectors $\ta,\tb$ and $f_{\kappa,V}(\ta)=f_{\kappa,V}(\tb)$ for some power-law kinetics $(\kappa,V)\in\mmK^{g}_{m,n}(Z)$. However, the latter  property could be fulfilled by many classes of kinetics other  than the class of power-law kinetics.

One such class of kinetics is \emph{Hill-type kinetics} that often is employed in modeling of biochemical reaction networks. In this section we will show that being $Z$-injective over $\mmK^g_{m,n}(Z)$ is equivalent to being injective over a  similar class of Hill-type kinetics. In particular, this implies that injectivity over $\mmK_{m,n}(Z)$ can be settled by applying Hill-type kinetics rather than power-law kinetics.

 We say that a kinetics is of Hill-type with respect to an $m\times n$ influence specification $Z$ if 
 $K=(K_1,\ldots,K_m)$  takes the form
$$K_j(c)=k_j \prod_{ i=1 }^n \frac{c_i^{v_{j,i}}}{\delta_{j,i}+c_i^{v_{j,i}}},$$
with $c\in\omR^n_+$ (defined by continuity at the boundary),  $k_j\in \R_+$, $\delta_j\in\omR^n_+$ and  $v_j\in\R^n$  for $j=1,\ldots,m$, such that 
$${\supp}^+(v_j)=z^+_j,  \quad  {\supp}^-(v_j)= z^-_j,\quad \text{and}\quad \supp(\delta_j)= \supp(v_j).$$
 The definition is very similar to that of power-law kinetics with the only difference being the factor $\delta_{j,i}$ in the denominators. A term with $v_{j,i}>0$ defines a positive influence,  
 while a term with $v_{j,i}<0$ defines a negative influence. Compared to power-law kinetics the  constant $\delta_{j,i}$ moderates a negative influence for low concentrations.

 Let $\kappa=(k_1,\ldots,k_m)$, $\mathbf{d}=(\delta_1,\ldots,\delta_m)$,  and  $V$ be an $m\times n$ matrix. We  denote a Hill-type kinetics by  $K=(\kappa,\mathbf{d},V)$, the set of Hill-type kinetics   with respect to  $Z$ by $\mmK^H_{m,n}(Z)$. 
Hill-type kinetics include Michaelis-Menten kinetics as a special case when  $v_{j,i}$ is one \cite{enz-kinetics}. 
 In contrast, power-law kinetics are not of Hill-type. However, power-law kinetics can be obtained as a limiting case of Hill-type kinetics by letting $k_{j}$ and the non-zero entries of $v_{j}$ tend to infinity such that $k_{j}/\prod_i v_{j,i}$ converges to a positive constant.
 
 Hill-type kinetics might be considered biochemically more reasonable than power-law kinetics as they are defined for all $\omR^n_+$ in contrast to  power-law kinetics that might not be defined for points on the boundary of $\omR^n_+$. In addition, Hill-type kinetics or Michaelis-Menten kinetics are often obtained when variables (species) are eliminated from the modelling equations \cite{enz-kinetics}.

\begin{theorem}\label{hill} 
Let $A$ be an $n\times m$ stoichiometric matrix, $Z$ and $m\times n$ influence specification and  $a,b\in\R^n_+$. Then:
\begin{enumerate}[(i)]
\item For every Hill-type kinetics  $K=(\kappa,\mathbf{d},V)\in\mmK^H_{m,n}(Z)$ there exists a
power-law kinetics $(\lambda,W)\in\mmK^{g}_{m,n}(Z)$ such that $Z(W)=Z$, 
$f_K(a)=f_{\lambda,W}(a)$ and $f_K(b)=f_{\lambda,W}(b)$.
\item For every power-law kinetics $(\lambda,W)\in\mmK^{g}_{m,n}(Z)$  there exists a
Hill-type kinetics  $K=(\kappa,\mathbf{d},V)\in\mmK^H_{m,n}(Z)$ such that $Z(V)=Z$,
$f_K(a)=f_{\lambda,V}(a)$ and $f_K(b)=f_{\lambda,V}(b)$.
\end{enumerate}
In particular, $f_{\lambda,W}(a)=f_{\lambda,W}(b)$ for some power-law kinetics $(\lambda,W)\in\mmK^{g}_{m,n}(Z)$ if and only if  $f_K(a)=f_K(b)$ for some Hill-type kinetics  $K=(\kappa,\mathbf{d},V)\in\mmK^H_{m,n}(Z)$ such that $Z(V)=Z(W)$.
\end{theorem}

As a consequence, injectivity of a stoichiomtric matrix $A$ over $\mmK_{m,n}(Z)$ is guaranteed by injectivity of $A$ over $\mmK^H_{m,n}(Z)$.   It also follows that $A$ is injective over $\mmK^g_{m,n}(Z)$  if and only if $A$ is injective over $\mmK^H_{m,n}(Z)$. Furthermore, we have that $A$ has multiple positive steady states in some stoichiometric class with respect to a Hill-type kinetics if and only if $A$ has multiple positive steady states in the same stoichiometric class with respect to a power-law kinetics.

\appendix

\section{Proofs}

\begin{proof}[Proposition~\ref{kernel}]  
We have that $v\in \ker(M)$ if and only if $Mv=0$. By assumption the top $d$ rows of $M$ are expressed as linear combinations of the bottom $s=n-d$ rows of $M$:
$$(0,\ldots, 0)=\omega^j M=M_{j}+\sum_{i=d+1}^n \omega^j_i M_{i},$$
where $M_{i}$ is the $i$-th row in $M$.
Hence, $v\in \ker(M)$ is equivalent to require that the scalar product of the bottom $s$ rows of $M$ and $v$ is zero. 
On the other hand, we have that  $v\in F$ if and only if $\omega^i\cdot v=0$ for all $i$.
It follows that $\ker(M)\cap F =\ker(\widetilde{M})$ and hence  $ \ker(M)\cap F =\ker(\widetilde{M})= \{0\} $ if and only if $\widetilde{M}$ has maximal rank, that is, if and only if $\det(\widetilde{M})\neq 0$. 
\end{proof}

\begin{proof}[Proposition~\ref{decomp}] 
Consider the characteristic polynomial of $M$ given as the determinant of $M-\lambda I_{n\times n}$. The right-hand side of the equality in  the proposition is $(-1)^d$ times the coefficient of $\lambda^{d}$. 
Let $P$ be the matrix whose top $d$ rows are $\omega^1,\dots,\omega^d$ and that agrees with the identity matrix in the bottom $s$ rows. By assumption $\omega^1,\dots,\omega^d$ is a reduced basis, hence the determinant of $P$ is $1$. Therefore:
$$\det(M-\lambda I_{n\times n}) = \det(PM-\lambda PI_{n\times n})=\det(PM-\lambda P). $$
Since the vectors  $\omega^1,\dots,\omega^d$ are orthogonal to the columns of $M$, the matrix $PM$ has zero rows in the top $d$ rows and agrees with $M$ in the bottom $s$ rows. It follows that 
$$PM-\lambda P =  \left( \begin{array}{c}-\lambda \omega^1 \\ \vdots \\ -\lambda\omega^d \\ M_{\{d+1,\dots,n\},\{1,\dots,n\}}-\lambda I_{s\times s} \end{array} \right).$$
The coefficient of $(-1)^d\lambda^d$ of the characteristic polynomial of $M$ is thus precisely given as the determinant of $\wm$.
\end{proof}

\begin{proof}[Proposition~\ref{coefs2}] 
Using \eqref{cauchy}, we have that
$$\det(J_c(\widetilde{f}_{\kappa,V})) = \sum_{I,J\subseteq \{1,\dots,n\}, \#I=\#J=s} \det(A_{I,J})\det((\partial K)_{J,I}).$$
We have, $(\partial K)_{j,i} = (k_{j}  c^{v_j}) c_i^{-1}v_{j,i}.$
That is, each term in  the $j$-th row of $\partial K$ is multiplied by $k_{j}  c^{v_j}$ and each term in the $i$-th column of $\partial K$ is multiplied by $c_i^{-1}$. 
It follows that
$$ \partial K = \diag(u_1,\dots,u_m) V \diag(g_1,\dots,g_n)$$
with $u_j=k_{j}  c^{v_j}$ and $g_i=c_i^{-1}$. Hence,
$$\det((\partial K)_{J,I}) = \prod_{j\in J}   u_j  \prod_{i\in I} g_i \det(V_{J,I}).    $$
The coefficient of  $\prod_{j\in J} k_{j}$ is given by summing these terms over all possible sets $I\subseteq\{1,\dots,n\}$
of cardinality $s$ and we obtain the expression in the statement.
\end{proof}

\begin{proof}[Proposition~\ref{det-criterion}]   
If (ii) holds then (i) is a consequence of  Theorem~\ref{injecclose-pl} and Proposition~\ref{coefs2}. To show that (i) implies (ii) consider the proof of Proposition~\ref{decomp} and the notation introduced there.
We have that 
$$\det(J_c(\widetilde{f}_{\kappa,V})) = \sum_{I,J\subseteq \{1,\dots,n\}, \#I=\#J=s} \det(A_{I,J})\det(V_{J,I})\prod_{j\in J}   u_j  \prod_{i\in I} g_i .$$
This is a linear polynomial in $u_*$ and $g_*$ and the coefficient of each monomial is given by the product of the two determinant. If two coefficients have opposite signs, then we can find values for $u_*$ and $g_*$ that make the polynomial vanish (Section 2). 
\end{proof}

\begin{proof}[Lemma~\ref{lemma-mono0}]
We first  prove the forward implication, that is,  that a strictly monotonic kinetics fulfills (i) and (ii).
Assume that  $K_{j}(a)>K_{j}(b)$. We show that:
$$(*) \quad \exists i: \quad a_i>b_i\hspace{.3cm} \text{ and } \hspace{.3cm} z_{j,i}=+, \quad \text{ or } \quad a_i<b_i \hspace{.3cm}\text{ and } \hspace{.3cm} z_{j,i}=-.$$
 To prove this, we assume that the contrary holds,  that is, for all $i$,
\begin{equation}\label{a2}
a_i\leq b_i, \quad \text{if} \quad z_{j,i}=+\qquad \textrm{and}\qquad a_i\geq b_i, \quad \text{if} \quad z_{j,i}=-.
\end{equation}
Let $\tilde{a}=a\wedge b$ be the minimum of $a$ and $b$,  $\tilde{a}_i=\min(a_i,b_i)$. By definition, $\tilde{a}\in\Omega_K$.
Recall that $K_{j}(c)$ does not depend on $c_i$ for which $z_{j,i}=0$ and $K_{j}(c)=0$ whenever $c\not\in\Omega_K(z^+_{j})$. Since $K$ respects the influence specification we have by monotonicity and \eqref{a2},
\begin{equation}\label{KK}
K_{j}(a)\leq K_{j}(\tilde{a})\leq  K_{j}(b).
\end{equation}
  However, this contradicts that 
$K_{j}(a)>K_{j}(b)$, implying that $(*)$ is true, hence also (i).

Assume now $K_{j}(a)=K_{j}(b)$. If $K_{j}(a)=0$ then there are $i,j\in z^+_{j}$  such that  $a_i=b_j=0$. Since  $a$ and $b$ are non-overlapping this cannot be the case and consequently  $K_{j}(a)\not=0$. It follows that   either  $a_i=b_i$ for all $i$ such that $z_{j,i}\not=0$, or  $a_i\neq b_i$ for some $i$ such that $z_{j,i}\not=0$.
If for all such indices  $\sign(a_i-b_i)z_{j,i}$ takes the same value $\epsilon=+$ or $-$, then by monotonicity  $\sign(K_{j}(a)-K_{j}(b))=\epsilon\neq 0$, which is a contradiction. Therefore, there exists two indices $i,\ell$  fulfilling (ii). It completes the first part of the proof.

To prove the reverse implication, assume that (i) and (ii) are fulfilled. Let $c,d\in \Omega_K(z_{j}^+)$ be two vectors that differ only in the $i$-th coordinate. If  $K_{j}(c)> K_{j}(d)$, then by $(i)$ we have $\sign(c_i-d_i)=z_{j,i}\neq 0$. It follows that if $z_{j,i}=0$, then $K_{j}(\cdot)$ is constant in the $i$-th coordinate. If $i\in  z^+_{j}\cup z^-_{j}$ and $K_{j}(c)= K_{j}(d)$ then according to (ii) we have $c_i=d_i$ contradicting $c_i\neq d_i$ (the second option cannot occur  since $c$ and $d$ differ in exactly one coordinate). Therefore 
$K_{j}(c)\neq K_{j}(d)$. Using (i) we conclude that $K_{j}(\cdot)$ is increasing/decreasing in the $i$-th coordinate depending on the sign of $z_{j,i}$. This completes the proof.
\end{proof}

\begin{proof}[Lemma~\ref{support}]   
(i) $Z$ has a signed $A$-determinant if and only if the polynomial $p_{Z}(X)$ has  constant sign when evaluated in positive values of the non-zero entries of $X$. The equivalence follows from the fact that each variable has degree zero or one in $p_{Z}(X)$.
(ii)  The polynomial $p_{Z'}(X)$ can be obtained from $p_{Z}(X)$ by setting some variables to zero. 
Statements (a)-(c) follow from this observation and statement (i).
\end{proof}

\begin{proof}[Proposition~\ref{openpos2}]  
By Lemma~\ref{subsets}, $A$ is injective over $\mmK^g_{m,n}(Z)$ if and only if $A$ is injective over $\mmK^g_{m,n}[V]$ for all $V$ such that  $Z(V)=Z$.
If $Z$ is A-SNS then each of the non-zero terms  $p_{Z,I,J}$ is sign-nonzero, have the same sign for all $I,J$ and at least one of them is non-zero. By Proposition~\ref{det-criterion} this implies that $A$ is injective over $\mmK^g_{m,n}[V]$ for all $V$.

Reciprocally, let us assume that $A$ is injective over $\mmK^g_{m,n}[V]$ for all $V$. Then by Proposition~\ref{det-criterion}, for each fixed $V$,   there is at least one non-zero term $p_{Z,I,J}(|V|)$, and all  nonzero terms have the same sign. It follows that $p_{Z}(|V|)\neq 0$ for all $V$. If there exists $V_1,V_2$  in  $\Sigma(Z)$ such that the sign of $p_{Z}(|V_1|)$ and $p_{Z}(|V_2|)$ are different, then by continuity we could find  $V_0$ such that  $p_{Z,I,J}(|V_0|)=0$  (the set of kinetic orders with associated influence $Z$ inherits a Euclidean topology from the Euclidean space it is embedded in). This contradicts that $A$ is injective over  $\mmK^g_{m,n}[V_0]$. Therefore, the sign of $p_{Z}(|V|)$ is independent of $V$ and thus by definition $Z$ is A-SNS.
\end{proof}

\begin{proof}[Proposition~\ref{openpos}] 
If $A$ is  injective over  $\bigcup_{Z \vert Z_1\preceq Z\preceq Z_2}\mmK^{g}_{m,n}(Z)$, then (iii) is trivially fulfilled and (ii) follows from Proposition~\ref{openpos2}. Proposition~\ref{openpos2} and Lemma~\ref{support}(ii) give that (iii) implies (i). Finally, if (ii) holds then by Lemma~\ref{support}(ii) we have that $Z_1$ has a signed A-determinant. It follows from Corollary~\ref{openpos3} that (iii) holds.
\end{proof}

\begin{proof}[Corollary~\ref{non-injective}]
Since the determinant $\det(J_c(\widetilde{f}_{\kappa,Y}))$ is not identically zero there is a term in its expansion in $\kappa$ with positive coefficient and a term with negative coefficient. Since, $Z_{\mmC}\preceq Z(V)$, 
 the terms in the polynomial expansion of $\det(J_c(\widetilde{f}_{\kappa,V}))$ in $\kappa$ cannot have all the same sign. Thus  
 Proposition~\ref{det-criterion} implies that $\mmK^g_{m,n}[V]$ is not injective.
\end{proof}

\begin{proof}[Theorem~\ref{subnet}]  
Assume that $A$ is injective over $\mmK^{g}_{m,n}[V]$. Proposition~\ref{det-criterion}(ii) implies that the non-zero products $\det(A_{I,J})\det(V_{J,I})$ have the same sign $\delta$ for all sets $I,J\subseteq\{1,\dots,n\}$ of cardinality $s$  and that at least one of the products is non-zero. 
For the matrix $A_{*,J}$ there is only one choice of column indices, namely the set $J=\{1,\dots,s\}$ ($s$ is the rank of $A_{*,J}$).  Observe that $\det(V_{J,I}) = \det(V_{\{1,\dots,s\},I})$.
If $\det(A_{I,J})\det(V_{J,I})=0$ for all $I$ then  all steady states (if there are any) are degenerate (equation \eqref{jaccomp}). If for some $I$ we have  $\det(A_{I,J})\det(V_{J,I})\neq 0$, then it follows from Proposition~\ref{det-criterion}   that $A_{*,J}$ is injective over $\mmK^{g}_{m,s}[V_{J,*}]$  and hence all steady states are non-degenerate.
\end{proof}

\begin{proof}[Theorem~\ref{equiv}]  
The proof is inspired by arguments presented in \cite{shinar-conc}. By Lemma~\ref{diff-to-weak}, (i)  implies (ii), and since power-law kinetics are differentiable with respect to the influence specification, (ii) implies (iii). (iii) and (iv) are equivalent according to Proposition~\ref{openpos2}. Let us prove that (iii) implies (i). 
 Assume that $A$ is injective over $\mmK^{g}_{m,n}(Z)$ but not $Z$-injective over $\mmK_{m,n}(Z)$. 
    Then there exists $K\in\mmK_{m,n}(Z)$ and distinct non-overlapping vectors $a,b\in\omR^n_+$ such that $\gamma:=a-b\in\im(A)$ and $f_K(a)=f_K(b)$. We seek a contradiction to the fact that $A$ is injective over $\mmK^{g}_{m,n}(Z)$, that is, we seek a power-law kinetics $(\kappa,V)$ such that $f_{\kappa,V}(\ta)=f_{\kappa,V}(\tb)$ for two vectors $\ta,\tb\in\R^n_+$ with $\ta-\tb\in\im(A)$.
    
 For  a positive constant $\delta>0$ and a positive constant vector $\zeta\in \R^n_+$, define
  $$\widetilde{K}_{j}(c)=K_{j}(c)+\delta, \quad \ta=a+\zeta, \quad \text{ and }\quad \tb=b+\zeta$$
such that $\widetilde{K}_{j}(a)-\widetilde{K}_{j}(b)=K_{j}(a)-K_{j}(b)$, $\ta-\tb=a-b$ and $\ta,\tb$ are positive vectors in $\R^n_+$. Therefore 
$\widetilde{K}_{j}(a)>\widetilde{K}_{j}(b)$ if and only if $K_{j}(a)>K_{j}(b)$, $a-b\in\im(A)$ if and only if $\ta-\tb\in\im(A)$, and $a_i>b_i$ if and only if $\ta_i>\tb_i$ (and similar for equality).
Since $a,b$ are non-overlapping,   $K_{j}(a)$ and $K_{j}(b)$ cannot both be zero for the same reaction.
Assume that we can find a kinetic order $V$ such that $Z(V)=Z$ and such that for all reactions
\begin{equation}\label{equality} 
\frac{\widetilde{K}_{j}(a)}{\widetilde{K}_{j}(b)}=\frac{ \ta^{v_{j}}}{\tb^{v_{j}}}=  \prod_{i=1}^{n} \left(\frac{\ta_i}{\tb_i}\right)^{v_{j,i}}.
\end{equation}
Then, if we define $\kappa$ by $k_{j}=\widetilde{K}_{j}(b)/\tb^{v_{j}}$, we have 
$$\widetilde{K}_{j}(b)=k_{j}\tb^{v_{j}}\qquad \textrm{and}\qquad \widetilde{K}_{j}(a) = \widetilde{K}_{j}(b)\frac{\ta^{v_{j}}}{\tb^{v_{j}}}= k_{j}\ta^{v_{j}},$$ 
and thus $f_K(a)=f_K(b)$ implies $f_{\kappa,V}(\ta)=f_{\kappa,V}(\tb)$.

Let us prove \eqref{equality}.
Assume that $K_{j}(a)>K_{j}(b)$. Then, by Lemma~\ref{lemma-mono0}, there exists $i$ for which $\sign(a_i-b_i)=z_{j,i}\not=0$.
Suppose that $a_i>b_i$ (that is, $\ta_i>\tb_i$) and $z_{j,i}=+$. Let $v_{j,\ell}= z_{j,\ell}\cdot \varepsilon$ for all $\ell\not=i$ and some positive $\varepsilon\in\R_+$. With this choice, \eqref{equality} holds if we can find $\varepsilon$ and $v_{j,i}>0$ (because $z_{j,i}=+$) such that
\begin{equation}\label{equality2}
1<\frac{\widetilde{K}_{j}(a)}{\widetilde{K}_{j}(b)}=\left(\frac{\ta_i}{\tb_i}\right)^{v_{j,i}}\prod_{\ell\not=i}\left(\frac{ \ta_{\ell}}{\tb_{\ell}}\right)^{ z_{j,\ell}\cdot \varepsilon}.
\end{equation}
Since the function $v\mapsto (\ta_i/\tb_i)^v$ is increasing ($\ta_i/\tb_i>1$), starts at $1$ and tends to infinity as $v$ increases, we can indeed find  $v_{j,i}>0$, potentially by choosing $\varepsilon$ small, such that \eqref{equality2} holds. The case $b_i<a_i$ is treated similarly. 

If  $K_{j}(b)>K_{j}(a)$ we proceed in the same way by interchanging the role of $a$ and $b$.
Finally, assume that $K_{j}(a)=K_{j}(b)$. Then, by Lemma~\ref{lemma-mono0}, either $a_i=b_i$ for all  $i\in z^+_{j}\cup z^-_{j}$, or  $\sign(a_i-b_i)=z_{j,i}\neq 0$ and $\sign(a_{\ell}-b_{\ell})=-z_{j,\ell}\not=0$ for some distinct $i,\ell$.
 In the first case the kinetic vector $v_{j}$ with $v_{j,i}=z_{j,i}\cdot 1$ fulfills equality \eqref{equality} and further $v_{j}$ satisfies 
 $Z(V)=Z$. In the second case, we have four scenarios depending on $z_{j,i}=+,- $ and $z_{j,\ell}=+,-$.  If $z_{j,i}=z_{j,\ell}=+$, then we can find, as above,
$v_{j,i},v_{j,\ell}>0$, such that
$$1=\frac{\widetilde{K}_{j}(a)}{\widetilde{K}_{j}(b)}=\left(\frac{\ta_i}{\tb_i}\right)^{v_{j,i}}\left(\frac{\ta_{\ell}}{\tb_{\ell}}\right)^{v_{j,\ell}}\prod_{u\not=i,j}\left(\frac{ \ta_u}{\tb_u}\right)^{ z_{j,u}\cdot \varepsilon},$$
because $\ta_i>\tb_i$ and $\ta_{\ell}<\tb_{\ell}$.
Hence, $v_{j}=(v_{j,1},\ldots,v_{j,n})$ with $v_{j,u}= z_{j,u}\cdot\varepsilon$, $u\not=i,\ell$, fulfills \eqref{equality}. The other three scenarios are treated in the same way. 
In conclusion, we can find a power-law kinetics such that $f_{\kappa,V}$ is not injective. 
\end{proof}

\begin{proof}[Theorem~\ref{injecclose}] 
Clearly (i) implies (ii).  Assume now that  (ii) holds. It is equivalent to $A$ being injective over $\mmK^g_{m,n}(Z)$, which again is equivalent to $Z$ being $A$-SNS (Proposition~\ref{openpos2}). Consider now condition (i). It is equivalent to $\det(J_c(\widetilde{f}_{K}))\not=0$ for all $c\in\R_+^n$ and  $K\in\mmK^d_{m,n}(Z)$ (equation \eqref{jaccomp}).  The Jacobian of $\widetilde{f}_K$ is $J_c(\widetilde{f}_{K}) =\widetilde{A(\partial K)}$,   where  $\partial K=\partial K(c)$ is the  $m\times n$ matrix with $\partial K_{j,i}= \partial K_{j}(c)/\partial c_i$. By definition of $K\in\mmK^d_{m,n}(Z)$, we have $z_{j,i}= \sign((\partial K)_{j,i})$. Since $Z$ is $A$-SNS, $\det(\widetilde{AZ(V)})\not=0$ for all kinetic orders $V$ with  $Z(V)=Z$. In particular this is true for the kinetic order given by $V=\partial K$. Hence condition (i) is true.
 \end{proof}

\begin{proof}[Proposition~\ref{graphical}]
Let $\mathfrak{G}_s$ denote the set of permutations of $s$ elements. By reordering  the species and reaction sets, we can assume that $I=J=\{1,\dots,s\}$. Then, by the definition of the determinant,
\begin{align*}
\det(A_{I,J})\det( (Z_X)_{J,I}) &= \sum_{\sigma,\tau \in \mathfrak{G}_s} \sign(\sigma)\sign(\tau) \prod_{i=1}^s a_{i,\sigma(i)} e_{\tau(i),i}.
\end{align*}
Fix a non-zero summand for some pair of permutations $\sigma,\tau$. Then $a_{i,\sigma(i)}\neq 0$ and $ e_{i,\tau(i)} \neq 0$ for all $i$. It follows that in $G_{A,Z}$ there is an edge from $r_{\sigma(i)}$ to $S_i$ and an edge from $S_{i}$ to reaction $r_{\tau(i)}$ for all $i$. Further, the set of these edges forms a $2s$-nucleus $D_{\sigma,\tau}$ with label  $ \prod_{i=1}^s a_{i,\sigma(i)} e_{\tau(i),i}$. Indeed,  each species node $S_i$ has precisely one ingoing edge with label $a_{i,\sigma(i)}$ and one outgoing edge with label $e_{\tau(i),i}$, and similarly for each reaction node $r_j$. 
Reciprocally, each $2s$-nucleus of $G_{A,Z}$ with vertices $S_1,\dots,S_s$, $r_1,\dots,r_s$ gives rise to a determinant term:
for each species node $S_i$ consider the ingoing and an outgoing edge $r_{j}\rightarrow S_i\rightarrow r_{j'}$ and define $\sigma(i)=j$ and $\tau(i)=j'$.

It remains to check that $\sign(\tau)\sign(\sigma)=\sign(D_{\sigma,\tau})$. The sign of $\tau\sigma$ agrees with the sign of $\tau\sigma^{-1}$, which in turn agrees with $(-1)$ to the number $p$ of cycles in the permutation  with even number of elements. Consider the graph in the reaction nodes $r_1,\dots,r_s$  obtained from $D_{\sigma,\tau}$ by removing the species nodes and joining two reaction nodes if they are connected through a species node. There is a correspondence between cycles of $\tau\sigma^{-1}$ and circuits in this collapsed graph. Therefore, the sign of  $\tau\sigma^{-1}$ is precisely $(-1)^p$.
\end{proof}

\begin{proof}[Lemma~\ref{gma-influence}]
 Assume that $K\in\mmK^w_{m,n}(Z)$ is a power-law kinetics with kinetic order $V$. Then $K\in\mmK^g_{m,n}(Z')$ with $Z'=Z(V)$. Let us prove that $Z'\preceq Z$. Let $i$ be such that $z_{j,i}'=+$, that is, $v_{j,i}>0$. Consider $a,b\in\R^n_+$ such that $a_u=b_u$, $i\not=u$, and $a_i>b_i$. Then  $K_{j}(a)>K_{j}(b)$ because $K$ is a power-law kinetics. By Definition~\ref{def-mono}(i) and using that $a,b$ only differ in the index $i$, we have $z_{j,i}=\sign(a_i-b_i)=+$. Therefore, $z_{j,i}'=z_{j,i}$. We proceed similarly if $z_{j,i}'=-$ to conclude that $Z'\preceq Z$.
\end{proof}

\begin{proof}[Theorem~\ref{equivw}]  
By Lemma~\ref{gma-influence},  (i) implies (iii). (ii) and (iii) are equivalent due to Theorem~\ref{equiv}. That (iii) implies (i) is proved similarly to the proof of Theorem~\ref{equiv}: The vector $v_{j}\in\R^n$ is likewise chosen
such that equation~\eqref{equality2} is fulfilled for the given influence specification $Z$.
\end{proof}

\begin{proof}[Theorem~\ref{soule}] 
The matrix $A$ is injective over $\mmK^d_{2s,n}(Z)$ if and only if $Z$ is $A$-SNS, that is, the non-zero coefficients of $p_Z(X)$ have constant sign and at least one is non-zero. The polynomial $p_Z(X)$ is the determinant of the symbolic matrix $\widetilde{AZ_X}$. Since the $n-s$ rows of $A$ are zero, a basis of $\im(A)^\perp$ is  $\{e^{s+1},\ldots,e^n\}$, where $e^j$ is the $j$-th unit vector (of length $n$). Then, an easy computation shows that  $p_Z(X)$ is the upper-left $s\times s$ minor of $AZ_X$. By construction, the sign pattern of the upper-left $s\times s$ minor of $AZ_X$ is  the upper-left  $s\times s$ submatrix of $\widehat{G}$, $\widehat{G}_L$.  

On the other hand, consider the interaction graph $G$. Only nodes $1,\dots,s$ have  incoming edges, hence a node $j>s$ cannot be part of any circuit of $G$. Consequently, any $s$-nucleus of $G$ contains precisely the nodes  $1,\dots,s$ and is a nucleus in the subgraph $G_L$ of $G$ given by these nodes. The matrix associated with this subgraph is $\widehat{G}_L$.

Let $t$ be  a non-zero term in the expansion of the upper-left $s\times s$ minor of $AZ_X$ and $N(t)$ the corresponding nucleus of  $G_L$. From \cite[Lemma 1]{soule} we have
$$\sign(N(t))=\sign(t)(-1)^{s+1}. $$
Consequently, all terms $t$ have the same sign if and only if all $s$-nuclei have the same sign, and there is a non-zero term if and only if there is an $s$-nucleus. Using Corollary~\ref{polyterms}, this proves the equivalence between (i) and (ii).
\end{proof}

\begin{proof}[Theorem~\ref{hill}]
Let $a,b\in\R^n_+$, $(\lambda,W)\in\mmK^{g}_{m,n}(Z)$ and $Z=Z(W)$, with $\lambda=(l_1,\dots,l_m)$.
For each $j=1,\ldots,m$, define $\delta_{j,i}=0$ and $v_{j,i}=0$ for all $i\in z_{j}^0$ (that is, for all $i$ such that $w_{j,i}=0$).  Let $M_{j}\in \R_+$  be such that $M_{j} a_i^{w_{j,i}}<1$ and $M_{j} b_i^{w_{j,i}}<1$ for  all $i\in  z_{j}^+\cup z_{j}^-$. Then we can find   $\delta_{j,i}\in\omR_+$ and $v_{j,i}\in\R$ such that,
\begin{equation}\label{hill-mass}
\frac{a_i^{v_{j,i}}}{\delta_{j,i}+a_i^{v_{j,i}}}=M_{j} a_i^{w_{j,i}}\quad \text{ and }\quad \frac{b_i^{v_{j,i}}}{\delta_{j,i}+b_i^{v_{j,i}}}=M_{j}b_i^{w_{j,i}}.
\end{equation}
Let $\mathbf{d}=(\delta_1,\ldots,\delta_m)$, $V=(v_{j,i})_{j=1,
\ldots,m,i=1,\ldots,n}$, and define $\kappa=(k_1,\ldots,k_m)$ by $k_{j}=l_{j} / M_{j}^J$ where $J$ is the cardinality of $z_{j}^+\cup z_{j}^-$. Then $Z(V)=Z(W)$
and $K=(\kappa,\mathbf{d},V)\in\mmK^H_{m,n}(Z)$. Further, 
$$l_{j} a^{w_{j}}= k_{j} \prod_{ i=1 }^n \frac{a_i^{v_{j,i}}}{\delta_{j,i}+a_i^{v_{j,i}} },$$
where $w_{j}=(w_{j,1},\ldots,w_{j,n})$, 
and similarly for $b$. It follows that $f_K(a)=f_{\kappa,W}(a)$ and $f_K(b)=f_{\kappa,W}(b)$.  This proves (ii). To prove (i) we follow the reverse procedure by choosing $M_{j}$ and $w_{j,i}$ to fulfill  equation~\eqref{hill-mass}. 
\end{proof}

\medskip
\textbf{Acknowledgements.}
EF has been supported by the postdoctoral grant ``Beatriu de Pin\'os'' from the Generalitat de Catalunya,  and project  MTM2012-38122-C03-01 from the Ministerio de Econom\'{\i}a y Competitividad of the Spanish government.  CW is supported by the Lundbeck Foundation, Denmark, The Danish Research Councils and the Leverhulme Trust, UK.  This work was initiated while EF and CW were visiting Imperial College London in fall 2011. The anonymous reviewers are thanked for their comments.

\end{document}